\documentclass[12pt]{amsart}
\usepackage[mathscr]{euscript}
\usepackage{amscd}
\newcommand{\nn}[1]{(\ref{#1})}

\newtheorem{theorem}{Theorem}[section]
\newtheorem{lemma}[theorem]{Lemma}
\newtheorem{proposition}[theorem]{Proposition}

\numberwithin{equation}{section}

\newcommand{\sfrac}[2]{{\textstyle \frac{#1}{#2}}}

\newcommand{\td}{\tilde{D}}
\newcommand{\nd}{\nabla}

\newcommand{\up}{\Upsilon}
\newcommand{\Up}{\up}
\newcommand{\ce}{\mathcal E}

\newcommand{\E}{\mathcal E}
\newcommand{\K}{\mathcal K}

\newcommand{\Co}{\mathcal C}

\renewcommand{\P}{\mathcal P}

\newcommand{\C}{\mathbb C}
\newcommand{\R}{\mathbb R}
\newcommand{\Z}{\mathbb Z}
\renewcommand{\Re}{\operatorname{Re}}
\renewcommand{\Im}{\operatorname{Im}}
\newcommand{\Ric}{\operatorname{Ric}}

\newcommand{\N}{\mathbb N}
\newcommand{\zb}{\overline z}

\newcommand{\Bb}{\overline{B}}

\newcommand{\Dt}{\widetilde \Delta}

\newcommand{\al}{\alpha}
\newcommand{\be}{\beta}
\newcommand{\beb}{\overline{\beta}}
\newcommand{\sib}{\overline{\sigma}}
\newcommand{\bb}{\overline{b}}

\newcommand{\La}{\Lambda}

\newcommand{\om}{\omega}
\newcommand{\g}{\gamma}
\newcommand{\gt}{\tilde{g}}
\newcommand{\s}{\sigma}

\newcommand{\bt}{{\boldsymbol \theta}}

\newcommand{\bh}{{\boldsymbol h}}

\newcommand{\pa}{\partial}

\newcommand{\spa}{\operatorname{span}}

\def\crn#1#2{{\vcenter{\vbox{
        \hbox{\kern#2pt \vrule width.#2pt height#1pt
           }
          \hrule height.#2pt}}}}
\def\intprod{\mathchoice\crn54\crn54\crn{3.75}3\crn{2.5}2}
\def\into{\mathbin{\intprod}}
\newsavebox{\ttoiso}
\sbox{\ttoiso}{\begin{picture}(30,12)(-5,0)
\put(0,3){$\simeq$}
\put(0,-3){$\longrightarrow$}
\end{picture}}

\newsavebox{\ccce}
\sbox{\ccce}{\begin{picture}(10,10)(0,0)
\put(0,0){$\ce$}
\put(2,0){$\ce$}
\end{picture}}

\newcommand{\bbox}                         %bbox=bbb box
{\mbox{$
\begin{picture}(9,8)(-.5,-1)
\put(0,6.4){\line(1,0){6.4}}
\put(0,0){\line(1,0){6.4}}
\put(0,0){\line(0,1){6.4}}
\put(1.5,0){\line(0,1){6.4}}
\put(6.4,0){\line(0,1){6.4}}
\end{picture}$}}

\newcommand{\bbbox}                         %bbox=bar over bbb box
{\ol{\bbox}}

\def\squarebox#1{\hbox to #1{\hfill\vbox to #1{\vfill}}}
\newcommand{\stopthm}{\hfill\hfill\vbox{\hrule\hbox{\vrule\squarebox
                 {.667em}\vrule}\hrule}\smallskip}

\newcommand{\ol}[1]{\overline{#1}}

\def\sideremark#1{\ifvmode\leavevmode\fi\vadjust{\vbox to0pt{\vss% the remark
 \hbox to 0pt{\hskip\hsize\hskip1em%                          will apper only
 \vbox{\hsize3cm\tiny\raggedright\pretolerance10000%          on the side
 \noindent #1\hfill}\hss}\vbox to8pt{\vfil}\vss}}}%
\def\thefootnote{}

\begin{document}

\title[Powers of the sub-Laplacian]{CR Invariant powers of the sub-Laplacian}
\author[A. R. Gover]{A. Rod Gover}
\address{Department of Mathematics\\
  University of Auckland\\
  Private Bag 92019\\
  Auckland, New Zealand}
 \email{gover@math.auckland.ac.nz}

\author[C. R. Graham]{C. Robin Graham}
\address{Department of Mathematics\\
University of Washington\\
Box 354350\\
Seattle, WA 98195 USA}
\email{robin@math.washington.edu}
  %\date{~/conf/formal.tex}

\begin{abstract}
CR invariant differential operators on densities with leading part a power
of the sub-Laplacian are derived.  One family of such operators is
constructed from the ``conformally invariant powers of the Laplacian'' via
the Fefferman metric; the powers which arise for these 
operators are bounded in terms of 
the dimension.  A second family is derived from a CR tractor calculus which
is developed here; this 
family includes operators for every positive power of the sub-Laplacian.
This result together with work of \u{C}ap, Slov\'{a}k and
Sou\u{c}ek imply in three dimensions 
the existence of a curved analogue of each such operator in flat space.
\end{abstract}

\renewcommand{\thefootnote}{}
\footnotetext{This research was partially supported by the Australian
Research Council and NSF grant DMS 0204480.  
This support and the hospitality of MSRI and the Universities of Auckland 
and Washington are gratefully acknowledged. }

\maketitle
\thispagestyle{empty}
\section{Introduction}

Invariant differential operators have a long history of importance
and this is particularly the case for operators of Laplace type.
The conformally invariant Laplacian is the basic example in conformal 
geometry.  A family of higher order generalizations of the
conformal Laplacian with principal part a power of the Laplacian 
was constructed in \cite{GJMS} and has been the subject of recent
interest.  In CR geometry, the CR invariant sub-Laplacian of Jerison-Lee
(\cite{JLee}) plays a role analogous to that of the 
conformal Laplacian.  In this paper we construct and study generalizations 
of the Jerison-Lee sub-Laplacian which are CR analogues of the 
``conformally invariant powers of the Laplacian''.  

One can deduce the existence of some such operators in the CR case from
the conformal operators via the Fefferman metric.  
Fefferman \cite{f} showed  that to a CR manifold $M$ of dimension $2n+1$
one can associate a conformal structure on a circle bundle $\Co$ of
dimension $N=2n+2$.
On conformal manifolds of dimension $N$, the construction in \cite{GJMS}
produces for each $w$, such that $N/2 + w =k\in \N$ and $k\leq N/2$
if $N$ is even, a 
conformally invariant natural differential operator
\mbox{$P_k:\ce(w)\rightarrow \ce(w-2k)$}, with principal part $\Delta^k$,
where $\ce(w)$ denotes the space of conformal densities of weight $w$.   
As we explain in the next section, on a CR manifold one can consider
CR densities $\ce(w,w')$ parametrized by \mbox{(bi-)weights} $(w,w') \in \C\times 
\C$. It is necessary that   
$w-w'\in \Z$ in order that $\ce(w,w')$ be well-defined and this is always
implicitly assumed.  
The space $\ce(w,w')$ of CR densities on $M$ can be naturally regarded as a 
subspace of the space of conformal densities $\ce(w+w')$ on $\Co$.
The operators $P_k$ for the Fefferman metric can be shown to preserve the
subspaces of CR 
densities and it is straightforward to calculate the principal part of the
induced operators.  One thereby obtains:
\begin{theorem}\label{ambient}
Let $n+1+w+w'=k\in \N$ with $k\leq n+1$.  There is a CR invariant
natural differential operator 
$$
P_{w,w'}:\ce(w,w')\rightarrow \ce(w-k,w'-k)
$$ 
whose principal part agrees with that of $\Delta_b^k$.
\end{theorem}
\noindent
In the case that $M$ is a 
flat CR manifold, the hypothesis $k\leq n+1$ is unnecessary in this
construction and in the flat case the operators $P_{w,w'}$ exist for all
$(w,w')$ such that $n+1+w+w'\in \N$.

In the conformal case, it is conjectured that if $N$ is even and $k>N/2$,
there is no natural operator with principal part $\Delta^k$ mapping
$\ce(-N/2+k)\rightarrow \ce(-N/2-k)$.  This has been established for $N=4$, 
$k=3$ (\cite{Grno}) and  $N=6$, $k=4$ (\cite{Wu2}).  Our main result is
the existence in the CR case of the following family of operators, 
which includes operators of higher orders.
Set $\N_0 = \N \cup \{0\}$.
\begin{theorem}\label{intromain}
For each $(w,w')$ such that $n+1+w+w'=k\in \N$ and
$(w,w') \notin \N_0\times \N_0$,
there is a CR invariant natural differential operator
$$
\P_{w,w'}: \ce(w,w')\to \ce(w-k,w'-k)
$$
whose principal part agrees with that of $\Delta_b^k$.
\end{theorem}
\noindent

In order to construct the operators $\P_{w,w'}$, we derive a tractor
calculus for CR geometry which we anticipate will be of independent
interest.  The main ingredients are to define a tractor
bundle, a complex vector bundle of rank $n+2$ over $M$, together with a CR
invariant connection, and an extension of this connection to a so-called
tractor $D$ operator between weighted tractor bundles.  We give a direct
derivation of these in terms of the pseudohermitian 
Tanaka-Webster connection induced by a choice of contact form on $M$,
parallel to the derivation in \cite{BEGo}   
in the projective and conformal cases.  The tractor bundle and connection 
can also be derived from the CR Cartan connection of Cartan, Tanaka,
Chern-Moser, but we found it preferable to proceed directly, especially 
for the tractor $D$ calculus needed for the construction of the
invariant operators.  See \cite{BEGo}, \cite{Esrni}, \cite{Gosrni}, 
\cite{TAMS} for further discussion and background about tractors.  Given
this machinery, the construction of 
the operators $\P_{w,w'}$ follows a similar construction of Eastwood and
Gover in the conformal case described in \cite{Gosrni}, involving iterating
tractor $D$'s.  The operators so constructed
are {\em strongly invariant} in the sense of \cite{Esrni},
i.e. they act invariantly not just on scalar densities but also on
density-valued tractors.  
Note that the tractor construction does not
yield an operator in the case $w=w'=0$ even though the construction via
the Fefferman metric does; this is the only weight for
which this occurs.  However, when $n=1$, a refinement of the
tractor construction does produce a CR invariant (but not necessarily
strongly invariant) operator
$\P_{0,0}:\ce(0,0)\rightarrow \ce(-2,-2)$, which by direct calculation can
be seen to agree with $P_{0,0}$.  Based on the situation in
the conformal case 
(\cite{GoPet}), we do not anticipate that the operators $P_{w,w'}$ and
$\P_{w,w'}$ agree in general when both are defined.  In fact, one must make
choices of orderings and conjugations in defining $\P_{w,w'}$, and
one expects that different choices will in general lead to different
operators.  In the flat case, the tractor construction (suitably
interpreted) also produces 
operators for all $(w,w')$ for which $n+1+w+w'\in \N$;
these do agree with the corresponding operators constructed via
the Fefferman metric.  

A construction of invariant natural operators is given in
\cite{CSSannals} 
for general parabolic geometries, of which CR geometry is a special case.  
This construction produces curved analogues of many invariant operators on
the flat homogeneous model, in the case when the flat operators 
are standard and act between irreducible bundles corresponding to integral 
non-singular highest weight vectors.  For CR geometry with $n\geq 2$, any
flat operator  
$P_{w,w'}$ satisfying these criteria has $\min \{w,w'\}\leq -n-1$.  
So in these dimensions Theorem~\ref{intromain} substantially extends the
results of \cite{CSSannals}.  However, the three dimensional case is
special:  the relevant parabolic is a Borel subalgebra, so all flat
operators are standard.  There do exist bundles $\ce(w,w')$ with $w+w'+2\in
\N$ corresponding to non-integral or singular highest weight vectors, for
which \cite{CSSannals} does 
not apply.  However, one can check that if $(w,w')\in \N_0\times \N_0$,  
then the highest weight vector for the bundle $\ce(w,w')$
is integral and non-singular.
We are not aware that it has been established in general that a given
specific operator constructed in \cite{CSSannals} has the desired principal
part, or even is nonzero.  However, it is possible to show
that this is the case for all the operators which arise in three
dimensional CR geometry. 
Therefore, in three dimensions Theorem~\ref{intromain} and \cite{CSSannals}
together yield the existence of a curved version of each of the flat
operators: 
\begin{theorem}
For each $(w,w')$ satisfying $w+w'+2=k\in \N$, there is a natural
differential operator on three dimensional
CR manifolds mapping $\ce(w,w')\to \ce(w-k,w'-k)$, whose principal part
agrees with that of $\Delta_b^k$.  
\end{theorem}
\noindent
Kengo Hirachi has informed us that he has established the existence of such
an operator in three dimensions mapping $\ce(1,1)\rightarrow \ce(-3,-3)$
via the CR ambient metric with ambiguity derived in \cite{hirachiannals}. 
We are grateful to Andi \u{C}ap for pointing out to us the special
applicability of \cite{CSSannals} in the three dimensional case.

In Section \ref{pseu} we review the basic facts about CR structures and
pseudohermitian geometry.  We derive the tractor bundle, its connection,
curvature, and the tractor $D$ operators in Section \ref{tractors}.  These
in turn are used in Section \ref{maint} to construct the operators
$\P_{w,w'}$.  The main point in proving Theorem \ref{intromain} is that
there are natural candidate operators which can be written down in terms of 
iterated tractor $D$ operators whose principal part is a multiple of
$\Delta_b^k$ and it must be determined when the multiple is nonzero.  This
leads to some numerology in calculating the multiple for flat
space.  Finally in Section \ref{amb} we review Lee's formulation of the
Fefferman metric and the GJMS construction of the conformal operators and
present the details of the proof of Theorem~\ref{ambient}.  We show how
this derivation via the Fefferman metric
can be reformulated in terms of Fefferman's ambient
metric and also in terms of the K\"ahler-Einstein metric of Cheng-Yau when 
$M$ is a hypersurface in $\C^{n+1}$, and
close with a brief discussion of $Q$-curvature in the CR case.  
Our presentation in this section in terms of the Fefferman metric
is influenced by the point of view explicated in
\cite{FH}. In Sections
\ref{maint} and \ref{amb} we also prove that the  
operators $\P_{w,w'}$ and $P_{w,w'}$ are self-adjoint (for $\P_{w,w'}$ the
orderings of the 
$D$ operators in the iteration must be chosen appropriately).  

It is a pleasure to thank Andi \u{C}ap, Mike Eastwood and Kengo Hirachi for
helpful discussions.  

\section{Pseudohermitian Geometry} \label{pseu}

A {\em CR-structure} on a $(2n+1)$-dimensional smooth manifold $M$
is an $n$-dimensional complex subbundle
$T^{1,0}\subset {\mathbb C}TM$ such that $T^{1,0}\cap \overline{T^{1,0}}
=\{0 \}$ and satisfying the integrability
condition $[T^{1,0},T^{1,0}]\subseteq T^{1,0}$, where we have used the same
notation $T^{1,0}$ for the bundle and the space of its smooth sections.
Set $T^{0,1}=\overline{T^{1,0}}$ and define
$\La^{1,0}\subset \C T^*M$ by $\La^{1,0} = (T^{0,1})^{\perp}$.
The {\em canonical bundle} $\K:=\La^{n+1}(\La^{1,0})$ is a complex line
bundle on $M$.  We will assume that $\K$ admits an $(n+2)^{\rm nd}$ root
and we fix a bundle denoted $\ce(1,0)$ which is a
$-1/(n+2)$ power of $\K$.  The bundle
$\ce(w,w'):=(\ce(1,0))^w\otimes (\ol{\ce(1,0)})^{w'} $
of $(w,w')$-densities is defined for $w,w'\in \C$ satisfying
$ w-w'\in {\mathbb Z}$.  (If $\ce(1,0)\setminus \{0\}$ is viewed as a
${\mathbb C}^*$-principal bundle, then $\ce(w,w')$ is the bundle induced by
the representation
$\lambda \rightarrow \lambda^{w}\ol{\lambda}{}^{w'}$ of
${\mathbb C}^*$.)  Without
further comment, we will assume implicitly whenever we write
$\E(w,w')$ that $w,w'\in \mathbb C$ and $w-w'\in\Z$.  
We will also usually write $\E(w,w')$ for
the space of sections of this bundle; the interpretation should be clear
from context.

The bundle
$$
H={\Re}(T^{1,0})
$$
is a real $2n$-dimensional subbundle of $TM$.  $H$ carries a natural almost
complex structure map given by
$J(Z+\bar{Z})=i(Z-\bar{Z})$ for $Z\in T^{1,0}$, which induces an
orientation on $H$.  We will assume that $M$ is orientable, which
implies that the bundle $H^{\perp}\subset T^*M$ admits a
nonvanishing global section.  A {\em pseudohermitian structure} on $M$ is
the choice of such a contact 1-form $\theta$.  We fix an orientation on the
bundle $H^{\perp}$ and restrict consideration to $\theta$'s which are
positive relative to this orientation.
The {\em Levi-form} of $\theta$ is the Hermitian form $h$ on $T^{1,0}$
defined by
$$
h(Z,{\ol{W}})=-2id\theta (Z,{\ol{W}}).
$$
We assume that the CR structure is nondegenerate, i.e. $h$ is a
nondegenerate form, whose signature we denote by $(p,q)$, $p+q=n$.
Given a pseudohermitian form $\theta$, define $T$ to
be the unique vector field on $M$ satisfying
\begin{equation}\label{Tnorm}
\theta (T)=1 ~~~{\rm and}~~~T\into d\theta=0.
\end{equation}
An {\em admissible coframe} is a set of $(1,0)$ forms $\{ \theta^{\al}\}$,
$\al =1,\cdots ,n$, which satisfy $\theta^{\al}(T)=0$ and whose
restrictions to $T^{1,0}$ form a basis for $(T^{1,0})^*$.
We will use lower case Greek indices to refer to frames for
$T^{1,0}$ or its dual.  We may also interpret these
indices abstractly, so will denote by
$\E^{\al}$ the bundle $T^{1,0}$ (or its space of sections) and by
$\E_{\al}$ its dual, and similarly for the conjugate bundles or for tensor
products thereof.  By integrability and \eqref{Tnorm}, we have
$$
d\theta = ih_{\al\beb} \theta^{\al}\wedge \theta^{\beb}
$$
for a smoothly varying Hermitian matrix $h_{\al\beb}$, which we may
interpret as the matrix of the Levi form in the frame $\theta^{\al}$, or
as the Levi form itself in abstract index notation.

It is shown in Lemma 3.2 of \cite{Lee1} (see also (4.3) of \cite{Lee2})
that if $\zeta$ is a (locally defined) nonvanishing section of
$\K$, then there is a unique pseudohermitian form $\theta$
with respect to which $\zeta$ is volume normalized in the sense that
$$
\theta\wedge(d\theta)^n = i^{n^2}n!(-1)^q\theta\wedge (T\into \zeta)
\wedge (T\into {\ol{\zeta}}).
$$
This is equivalent to the requirement that $\zeta$
have unit length with respect to the natural norm
induced by $\theta$ and $h$.
Upon replacing $\zeta$ by $\lambda \zeta$, where $\lambda$ is a smooth
$\C^*$-valued function, $\theta$ is replaced by
$|\lambda|^{2/{n+2}} \theta$.  Now
$|\zeta|^{-2/{n+2}}:=\zeta^{-1/(n+2)}\otimes \ol{\zeta}^{-1/(n+2)}$
is a section of
$\K^{-1/{n+2}}\otimes \ol{\K}^{-1/{n+2}}=\E(1,1)$, so
$\bt :=\theta \otimes |\zeta|^{-2/{n+2}}$ defines a section of
$T^*M \otimes \E(1,1)$ which depends only on the CR structure.
The Levi form of $\theta$ scales the same
way, so defines a canonical section ${\bh}_{\al\beb}$
of $\E_{\al\beb}(1,1)$, where in general we denote the tensor product with
a density bundle by omitting the second $\E$:
$\E_{\al}(w,w'):=\E_{\al} \otimes \E(w,w')$.  We will use exclusively
${\bh}_{\al\beb}$ and its inverse
${\bh}^{\al\beb}\in \E^{\al\beb}(-1,-1)$ to raise and lower indices, so
that raising and lowering indices changes weight.
The choice of pseudohermitian form
$\theta$ on $M$ is equivalent to the choice of a {\em CR-scale}
$0<t \in \E(1,1)$ related by $\bt = t\theta$, and we also have
${\bh}_{\al\beb}=t h_{\al\beb}$.  If $\Up\in C^{\infty}(M)$, the
change of pseudohermitian structure $\hat{\theta}=e^{\Up}\theta$
results in $\hat{h}_{\al\beb}=e^{\Up}h_{\al\beb}$ and
$\hat{t}=e^{-\Up}t$.
A pseudohermitian structure is said to be {\em flat} if it is locally
equivalent to the Heisenberg group with its standard pseudohermitian
structure.  It is said to be {\em CR flat} if there is a rescaling
of $\theta$ which is pseudohermitian flat.

A choice of pseudohermitian structure determines a connection on $TM$,
the Tanaka-Webster, or pseudohermitian, connection (\cite{tan},
\cite{web}).
It is given in terms of an admissible coframe by
\begin{equation}\label{connforms}
\nd \theta^{\al} = -\om_{\be}{}^\al \otimes \theta^{\be}, \quad \nd \theta
= 0,
\end{equation}
where Webster's connection 1-forms $\om_{\be}{}^{\al}$ satisfy
$$
d\theta^{\al}=\theta^{\be}\wedge\om_{\be}{}^{\al}+A^{\al}{}_{\beb}
\bt\wedge\theta^{\beb},
$$
and the {\em pseudohermitian torsion tensor}
$A_{\al\be}=\ol{A_{\ol{\al}\ol{\be}}}\in
\E_{(\al\be)}$ is symmetric.  (According to our conventions for raising
indices, $A^{\al}{}_{\beb}$ is a section of $\E^{\al}{}_{\beb}(-1,-1)$.)
In particular, this connection preserves $T^{1,0}$, so induces connections
on $\E^{\al}$ and $\E_{\al}$.  One has $\nd h = 0$.
There is an induced connection on the canonical bundle, and therefore
also on all the density bundles $\E(w,w')$.
\begin{proposition}\label{parallel}
$\nd \bt = 0$
\end{proposition}
\begin{proof}
The definition is $\bt = \theta \otimes |\zeta|^{-2/{n+2}}$, where
$\zeta$ is
volume normalized with respect to $\theta$.  Since $\nd \theta = 0$, it
suffices to show that $\nd |\zeta|^2=0$ for such $\zeta$.  Choose an
admissible coframe $\{\theta^{\al}\}$ such that
$\zeta = \theta\wedge\theta^1\wedge\cdots\wedge\theta^n$.  The
condition that $\zeta$ is volume normalized is equivalent to
$\det(h_{\al\beb})=(-1)^q$.  By  \eqref{connforms} we have
\begin{equation}\label{ndzeta}
\nd \zeta = -\om_{\al}{}^{\al} \otimes \zeta.
\end{equation}
It follows that
$\nd |\zeta|^2=-(\om_{\al}{}^{\al} + \ol{\om_{\al}{}^{\al}})\otimes
|\zeta|^2$.
However, from $\nd h = 0$ one obtains
$\om_{\al}{}^{\al} + \ol{\om_{\al}{}^{\al}}=h^{\al\beb}dh_{\al\beb}=0$ as
desired.
\end{proof}
\noindent
{From} Proposition~\ref{parallel} we conclude that $\nd t=0$, where $t$ is
the scale associated to $\theta$, and that $\nd \bh = 0$.

The pseudohermitian connection preserves the splitting
$
\C TM = T^{1,0}\oplus T^{0,1} \oplus \spa{T}.
$
Therefore, if we decompose
a tensor field relative to this splitting (and/or its dual), we may
calculate the covariant derivative componentwise.  Each of the components
may be regarded as a section of a tensor product of $\E^{\al}$ or its dual
or conjugates thereof.  Therefore we will often restrict consideration to
the action of the connection on $\E^{\al}$ or $\E_{\al}$.
We will use indices $\al, \ol{\al}, 0$ for components with
respect to the frame $\{\theta^{\al},\theta^{\ol{\al}}, \bt\}$ and its
dual, so that the 0-components incorporate weights.  If
$f$ is a (possibly density-valued)
tensor field, we will denote components of the (tensorial) iterated
covariant derivatives of $f$ in such a frame by preceding $\nd$'s, e.g.
$\nd_{\al} \nd_0 \cdots \nd_{\beb} f$.  As usual, such indices may
alternately be interpreted abstractly.  So, for example, if
$f_\be\in \E_\be(w,w')$, we will consider $\nd f$ as the triple
$\nd_{\al}f_\be\in \E_{\al\be}(w,w')$,
$\nd_{\ol{\al}}f_\be\in \E_{\ol{\al}\be}(w,w')$,
$\nd_0f_\be \in \E_\be(w-1,w'-1)$.

The torsion and curvature of $\nd$ can be described
by structure equations for the forms $\om_{\al}{}^{\be}$
(see \eqref{connforms} above and (1.3) of \cite{JLee2}),
or by commuting second
derivatives on functions and 1-forms (see Lemma 2.3 of \cite{Lee2}.
We generally follow the same conventions as \cite{Lee2},
\cite{JLee2}, except that we substitute $\bt, \bh$ for $\theta, h$ so
that all quantities are naturally weighted.)  These can be expressed in
terms of
the {\em pseudohermitian curvature tensor} $R_{\al\beb\g\ol{\delta}}\in
\E_{\al\beb\g\ol{\delta}}(1,1)$, $A_{\al\be}$,
$\nd_{\ol{\gamma}}A_{\al\be}$, and conjugates of the latter two.  Flat
pseudohermitian structures are characterized by the vanishing of
$R_{\al\beb\g\ol{\delta}}$ and $A_{\al\be}$.
The Webster-Ricci tensor is
defined by
$$
R_{\al\beb}=R_{\g}{}^{\g}{}_{\al\beb}\in \E_{\al\beb}
$$
and the Webster scalar curvature by
$$
R=R_{\al}{}^{\al}\in \E(-1,-1).
$$
{From} these we define
$$
P_{\al\beb}:=\frac{1}{n+2}\left(R_{\al\beb}
-\frac{1}{2(n+1)}R\bh_{\al\beb}\right).
$$
We will also need to know how to commute derivatives of densities.

\begin{proposition}\label{commprop}
If $f\in\E(w,w')$, then
\begin{equation}\label{dencomm}
\begin{split}
\nd_{\al} \nd_{\be} f -\nd_{\be} \nd_{\al} f &= \,\,0\\
\nd_{\al} \nd_{\beb} f -\nd_{\beb} \nd_{\al} f &
=\frac{w-w'}{n+2}R_{\al\beb}f-i\bh_{\al\beb}\nd_0 f\\
\nd_{\al} \nd_{0 } f -\nd_{0} \nd_{\al } f &=
\frac{w-w'}{n+2}(\nd^{\g}A_{\g\al})f+A_{\al\g}\nd^{\g}f.
\end{split}
\end{equation}
\end{proposition}
\begin{proof}
We first note that \eqref{dencomm} holds when $w=w'=0$, since in that case
it agrees with Lemma 2.3 of \cite{Lee2}.  Now it is a straightforward
consequence of the definitions that
if $V$ is a vector bundle with connection over a manifold $M$, which
itself has a linear connection, if $v^A$ is a section of $V$, then
\begin{equation}\label{gencomm}
[\nd_i,\nd_j]v^A = \Omega_{ijB}{}^A v^B - T_{ij}^k\nd_kv^A,
\end{equation}
where here $i,j,k$ label $TM$, $\Omega$ denotes the curvature
of the connection on $V$, $T_{ij}^k$ denotes the
torsion of the connection on $TM$, and the second
covariant derivative is with respect to the coupled connection on
$T^*M\otimes V$.
We first apply \eqref{gencomm} taking $V$ to be the trivial bundle with the
flat connection over our pseudohermitian manifold with its connection:
the case $w=w'=0$.  The curvature term vanishes, and we
conclude that the torsion term in \eqref{gencomm} is
exactly the right hand side of \eqref{dencomm} for $w=w'=0$.  Next
let $w=-(n+2), w'=0$, and take $V=\E(w,w')=\K$ with its
induced connection.  The form of the torsion term in \eqref{gencomm} is
independent of the
choice of $V$, so takes exactly the same form as in the previous case.
To identify the curvature term, observe that since $\K$ is a line bundle,
its curvature is a scalar 2-form
on $M$.  Choose an admissible coframe and set $\zeta =
\theta\wedge\theta^1\wedge\cdots\wedge\theta^n$.  Then we
have \eqref{ndzeta}, so the connection form for $\K$
relative to the frame $\{\zeta\}$ is $-\om_{\al}{}^{\al}$.  The
curvature of $\K$ is therefore $-2d\om_{\al}{}^{\al}$.  However,
(2.3), (2.4) of \cite{Lee2} state that
$$
d\om_{\al}{}^{\al}=R_{\al\beb}\theta^{\al}\wedge \theta^{\beb}
+\nd^{\be}A_{\beta\al}\theta^{\al}\wedge \theta
-\nd^{\beb}A_{\beb\ol{\al}}\theta^{\ol{\al}}\wedge \theta.
$$
Thus \eqref{gencomm} reduces to \eqref{dencomm} for $w=-(n+2), w'=0$.
The general case now follows upon conjugating, taking powers, and adding.
\end{proof}
\noindent
Proposition~\ref{commprop} together with the formulae in Lemma
2.3 of \cite{Lee2} enable one to calculate the effect of
commuting covariant derivatives on any weighted tensor field.

If $\hat{\theta}=e^{\Up}\theta$ is another pseudohermitian form, we can
express the connection, curvature, and torsion for $\hat{\theta}$ in terms
of that for $\theta$.  We consider first the relation between the induced
connections on tensor products of $\E^{\al}$ or its dual or
conjugates thereof.  For efficiency in expressing the
transformation laws, we will denote covariant derivatives of $\Up$
with indices:  $\nd_{\al}\nd_{\beb}\Up = \Up_{\beb\al}$.
One observes first that if $\{\theta^{\al}\}$ is an admissible coframe for
$\theta$, then
$\{\hat{\theta}^{\al}=\theta^{\al} +i\Up^{\al}\theta\}$ defines an
admissible coframe for $\hat{\theta}$.  Since the restrictions of
$\hat{\theta^{\al}}$ and $\theta^{\al}$ to $T^{1,0}$ agree, the
components of a section $f$ of $\E^{\al}$ or its dual (or conjugates or
tensor products) are the same in the two frames.  We will denote
by $\widehat{\nd}_{\al}f$, $\widehat{\nd}_{\ol{\al}}f$, and
$\widehat{\nd}_0 f$ the components of
$\widehat{\nd} f$ relative to
$\{\hat{\theta^{\al}},\hat{\theta^{\ol{\al}}}, \bt\}$.
For example, if $f$ is an unweighted function, one obtains upon changing
frames that
\begin{equation}\label{function}
\widehat{\nd}_{\al}f=\nd_{\al}f,\quad
\widehat{\nd}_{\ol{\al}}f=\nd_{\ol{\al}}f,\quad
\widehat{\nd}_0 f=\nd_0 f
+i\Up^{\ol{\g}}\nd_{\ol{\g}}f-i\Up^{\g}\nd_{\g}f.
\end{equation}
The connection forms $\hat{\om}_{\al}{}^{\be}$ for $\widehat{\nd}$
can be expressed in terms of $\om_{\al}{}^{\be}$ and
$\Up$; this is Lemma 3.4 of \cite{JLee2}.  It is straightforward to
use this to calculate the analogue of \eqref{function} for a section
$\tau_{\be}$ of $\E_{\be}$.  One obtains:
\begin{equation}\label{vector}
\begin{array}{rcl}
\widehat{\nd}_{\al} \tau_{\be}&=& \nd_{\al}
\tau_{\be}-\up_{\be}\tau_{\al}-\up_{\al}\tau_{\be} \\
\widehat{\nd}_{\bar{\al}} \tau_{\be}&=& \nd_{\bar{\al}} \tau_{\be}
+\bh_{\be\bar{\al}} \up^{\g} \tau_{\g}\\
\widehat{\nd}_0 \tau_{\be}&=&
\nd_0 \tau_{\be}+i\up^{\bar{\g}}\nd_{\bar{\g}} \tau_{\be}
-i\up^{\g}\nd_{\g} \tau_{\be}
-i(\up^{\g}{}_{\be} -\up^{\g}\up_{\be})\tau_{\g}.
\end{array}
\end{equation}
We also need to know how the connection transforms on densities.
\begin{proposition}\label{denhat}
If $f\in \E(w,w')$, then
\[
\begin{split}
\widehat{\nd}_{\al}f =&\nd_{\al}f + w\Up_{\al}f\\
\widehat{\nd}_{\ol{\al}}f = &\nd_{\ol{\al}}f +w'\Up_{\ol{\al}} f\\
\widehat{\nd}_0 f = &\nd_0 f
+i\Up^{\ol{\g}}\nd_{\ol{\g}}f-i\Up^{\g}\nd_{\g}f\\
&+\sfrac1{n+2}\left[(w+w')\Up_0+iw\Up^{\g}{}_{\g}
-iw'\Up^{\ol{\g}}{}_{\ol{\g}} +i(w'-w)\Up^{\g}\Up_{\g}\right]f
\end{split}
\]
\end{proposition}
\begin{proof}
As in the proof of Proposition~\ref{commprop}, if we establish the
result for
$w=-(n+2),w'=0$, the general case follows by conjugating, taking powers,
and adding.  Let $\zeta = \theta\wedge\theta^1\wedge\cdots\wedge\theta^n$
be a section of $\K=\E(-(n+2),0)$.
Then $\nd\zeta$ is given by \eqref{ndzeta}.  However, also
$\zeta =
e^{-\Up}\hat{\theta}\wedge\hat{\theta^1}\wedge
\cdots\wedge\hat{\theta^n}$, so applying \eqref{ndzeta} again gives
$\widehat{\nd}\zeta = -(\hat{\om}_{\al}{}^{\al} + d\Up)\otimes \zeta$.
Therefore
$$
\widehat{\nd}\zeta - \nd \zeta = -(\hat{\om}_{\al}{}^{\al} -
\om_{\al}{}^{\al} +d\Up)\otimes \zeta.
$$
Now Lemma 3.4 of \cite{JLee2} gives $\hat{\om}_{\al}{}^{\be} -
\om_{\al}{}^{\be}$ in terms of $\Up$.  Contracting, adding
$d\Up$, and reformulating the result in terms of components yields
the desired formulae.
\end{proof}

The curvature and torsion also transform under the pseudohermitian change.
{From} Lemma 2.4 of \cite{Lee2} one obtains
\begin{equation}\label{Ptrans}
\widehat{P}_{\al\beb}=P_{\al\beb}-
\sfrac{1}{2}(\up_{\al\beb}+\up_{\beb\al})
-\sfrac{1}{2}\up_{\g}\up^{\g} \bh_{\al\beb}
\end{equation}
and
\begin{equation}\label{Atrans}
\widehat{A}_{\al\be}=A_{\al\be}+i\up_{\al\be}
-i\up_{\al}\up_{\be}.
\end{equation}
We also will need the transformation laws for two other objects.
Set $P=P_{\al}{}^{\al}=\sfrac{1}{2(n+1)}R$ and define
$$
T_\al =\frac1{n+2}(\nd_\al P -i\nd^\be A_{\al\be})\in \E_{\al}(-1,-1),
$$
$$
S=-\frac{1}{n}(\nd^\al T_\al+\nd^{\bar{\al}}T_{\bar{\al}}
+P_{\al\beb}P^{\al\beb}-A_{\al\be}A^{\al\be})\in \E(-2,-2).
$$
Straightforward calculation using the formulae discussed above
and the Bianchi identities of Lemma 2.2 of \cite{Lee2} gives:
\[
\begin{array}{ll}
\begin{split}
\widehat{T}_\al = T_\al &+ \sfrac{i}2 \Up_{0\al}+P_\al{}^{\be}\Up_\be
-iA_{\al\be}\Up^\be \\& +\sfrac12\Up_{\al\be}\Up^\be
-\sfrac12\Up_\al{}^\be \Up_\be -\sfrac12\Up_\be\Up^\be\Up_\al,
\end{split}\\
\begin{split}
\widehat{S}= S &
+ \sfrac12\Up_{00}-3(\Up^\al T_\al
+\Up^{\ol{\al}}T_{\ol{\al}}) \\&
+i(\Up_{0\ol{\al}}\Up^{\ol{\al}} - \Up_{0\al}\Up^\al)
-\sfrac14(\Up_0)^2
+\sfrac{3i}2(A_{\al\be}\Up^\al\Up^\be
- A_{\ol{\al}\ol{\be}}\Up^{\ol{\al}}\Up^{\ol{\be}}) \\&
-3P_{\al\ol{\be}}\Up^\al\Up^{\ol{\be}}
-\sfrac12(\Up_{\al\be}\Up^\al\Up^\be +
\Up_{\ol{\al}\ol{\be}}\Up^{\ol{\al}}\Up^{\ol{\be}})\\&
+\sfrac12(\Up_{\al\ol{\be}}+\Up_{\ol{\be}\al})\Up^\al\Up^{\ol{\be}}
+\sfrac34(\Up_\al\Up^\al)^2.
\end{split}
\end{array}
\]
We derived the expressions for $T_\al$ and $S$ from the condition that they
satisfy transformation laws of this form, which we needed to construct the
tractor connection in the next section.  Only later did we observe that
they occur already in \cite{Lee1} as components of the 
connection forms and Ricci tensor of the Fefferman metric.
Also, $T_{\al}$ and a variant of $S$ arise in \cite{JLee2}, 
where they were combined with $P_{\al\beb}$ and
$A_{\al\be}$ to form a 2-tensor on $M$ which was used to construct a
version of pseudohermitian normal coordinates.  In that work, only the
part of the transformation laws involving highest derivatives
(counted non-isotropically) of
$\Up$ was relevant.  The version of $S$ in \cite{JLee2} therefore does not
include the $P_{\al\beb}P^{\al\beb}-A_{\al\be}A^{\al\be}$ term, which
however is important for us.

\section{Tractors} \label{tractors}

The pseudohermitian connection depends on
the choice of $\theta$.  In this section we construct a vector
bundle of rank $n+2$ over $M$, the tractor bundle, together with a CR
invariant tractor connection.  We also show how to extend this connection
to a CR invariant tractor $D$ operator on weighted tractors which can be
iterated.

For a given choice of $\theta$, the (co-)tractor bundle $\ce_{A}$ is
realized as a direct sum
\begin{equation}\label{dirsum}
\ce_A = \ce(1,0)\oplus\ce_\al(1,0)\oplus\ce(0,-1).
\end{equation}
The realization corresponding to $\hat{\theta}=e^\Up \theta$ is identified
with that for $\theta$ via
\begin{equation}\label{identif}
     \left( \begin{array}{c} \hat{\sigma}\\
                            \hat{\tau_{\al}}\\
          \hat{\rho} \end{array}\right)
=
\left( \begin{array}{c} \sigma\\
                            \tau_{\al}+\up_\al\sigma\\
\rho -\up^\be\tau_\be-\frac{1}{2}(\up^\be\up_\be+i\up_0)\sigma
\end{array}\right)
=
M_\Up  \left( \begin{array}{c} \sigma\\
                            \tau_{\be}\\
                 \rho \end{array}\right),
\end{equation}
where $\s\in \ce(1,0)$, $\tau_\be \in\ce_\be(1,0)$, $\rho\in\ce(0,-1)$,
and
$$
M_\Up=
\left(\begin{matrix}
1&0&0\\
\Up_\al&\delta_\al{}^\be&0\\
-\sfrac12(\Up^\g\Up_\g + i\Up_0)&-\Up^\be&1
\end{matrix}\right).
$$
Recalling \eqref{function}, it is easily checked that these identifications
are consistent upon changing to yet a third $\theta$.
We may therefore mod them out to
obtain a bundle $\E_A$ determined solely by the CR structure.
(More formally, the total space of the bundle $\E_A$ can be defined to be
the disjoint union of one copy of the total space of
$\ce(1,0)\oplus\ce_\al(1,0)\oplus\ce(0,-1)$ for each global section
$\theta$, modulo the equivalence relation \eqref{identif}.  The smooth
structure on this total space and the linear structure on the fibers are
inheirited from that of any of the realizations.)
The subspaces $\{\sigma = 0\}$ and $\{\sigma=\tau_\al=0\}$ are preserved
by the identifications, so determine canonical subbundles of $\E_A$.
This composition series structure is what remains of the direct sum
decomposition \eqref{dirsum} after making the identifications.
The vector given by $\sigma=\tau_\al=0,\rho=1$ is preserved by $M_\Up$, so
defines an invariant section of $\E_A(0,1)=\E_A\otimes \E(0,1)$
which we denote by $Z_A$.  The matrix $M_\Up$ is in $SU(h_{A\ol{B}})$,
where
$$
h_{A\ol{B}}=
\left(\begin{matrix}
0&0&1\\
0&\bh_{\al\beb}&0\\
1&0&0
\end{matrix}\right).
$$
It follows that $h_{A\ol{B}}$ defines an invariant Hermitian metric on
$\E_A$, and that $\E_A$ has an invariant volume form.
We denote by $\E_{\ol{A}}$ the conjugate bundle, by $\E^A$ the dual bundle
(the tractor bundle), and we use $h_{A\ol{B}}$
and its inverse to raise and lower tractor indices.

We define a connection on $\E_A$ as follows.  If
$$
v_A=\left( \begin{array}{c} \sigma\\
                            \tau_\al\\
                             \rho \end{array}\right)
$$
is a section of $\E_A$ in the realization determined by $\theta$, set
$$
\nd_\be v_A=\left( \begin{array}{c}
                            \nd_\be\sigma-\tau_\be\\
                            \nd_\be\tau_\al+iA_{\al\be}\sigma\\
         \nd_\be \rho-P_\be{}^\al\tau_\al+T_\be\sigma
                                   \end{array}\right) ,
$$

\begin{equation}\label{tracconn}
\nd_{\beb} v_A= \left( \begin{array}{c}
                               \nd_{\beb} \sigma\\
\nd_{\beb}\tau_\al
+\bh_{\al\beb}\rho+P_{\al\beb} \sigma\\
\nd_{\beb}\rho+iA_{\beb}{}^\al \tau_\al-T_{\beb}\sigma
 \end{array} \right),
\end{equation}
$$
\nd_0 v_A= \left( \begin{array}{c}
                               \nd_0 \sigma+\frac{i}{n+2}P\sigma-i\rho\\
\nd_0\tau_\al-iP_\al{}^\be\tau_\be+\frac{i}{n+2}P\tau_\al+
2iT_\al\sigma\\
\nd_0\rho+\frac{i}{n+2}P\rho+2iT^\al\tau_\al+iS\sigma
 \end{array} \right),
$$
where $T_\al$ and $S$ were defined in \S 2, and where the $\nd$'s on the
right hand side refer to the pseudohermitian connection on the appropriate
weighted bundles.  This operation satisfies a Leibnitz formula so defines a
connection on $\E_A$.  It is a direct calculation
using the transformation formulae in the previous section to show that
the definition is independent of the choice of $\theta$, and therefore CR
invariant.  For example, consider the second component of
$\nd_\be v_A$.  Using \eqref{vector}, Proposition~\ref{denhat}, and
\eqref{Atrans}, we obtain
\[
\begin{split}
\widehat{\nd}_\be  \hat{\tau}_\al + i\hat{A}_{\al\be}\hat{\s}
= &\widehat{\nd}_\be(\tau_\al+\Up_\al\s )+i\hat{A}_{\al\be}\s\\
=&\nd_\be(\tau_\al+\Up_\al\s)-\Up_\be(\tau_\al+\Up_\al\s)
-\Up_\al(\tau_\be+\Up_\be\s)\\
&\quad + \Up_\be(\tau_\al +\Up_\al\s)
+i(A_{\al\be}+i\Up_{\al\be}-i\Up_\al\Up_\be)\s\\
=&(\nd_\be\tau_\al + iA_{\al\be}\s)+\Up_\al(\nd_\be\s-\tau_\be)
\end{split}
\]
as desired.  (Here we verified the full transformation law, but it suffices
to check invariance to first order in
$\Up$.)  This connection induces connections on the dual and conjugate and
tensor product bundles, and a calculation shows that $\nd h_{A\ol{B}}=0$
(reflecting a Hermitian symmetry inherent in \eqref{tracconn}).
Therefore differentiation commutes with raising and lower tractor indices.

One can calculate the curvature of the tractor connection from
\eqref{gencomm}, using the commutation formulae and Bianchi identities
of \cite{Lee2}.  The result is the following:
$$
\Omega_{\rho\sigma A}{}^B=0, \qquad  \Omega_{{\ol{\rho}}\sib A}{}^B=0,
\qquad
\Omega_{\rho\sib A}{}^B=
\left(\begin{matrix}
0&0&0\\
iV_{\rho\sib\al}&S_{\rho\sib\al}{}^{\beta}&0\\
U_{\rho\sib}&-iV_{\sib\rho}{}^{\beta}&0
\end{matrix}\right),
$$
$$
\Omega_{\rho 0 A}{}^B=
\left(\begin{matrix}
0&0&0\\
Q_{\al\rho}&V_{\rho}{}^{\beta}{}_{\al}&0\\
Y_{\rho}&-iU_{\rho}{}^{\beta}&0
\end{matrix}\right),
\qquad
\Omega_{{\ol{\rho}}0 A}{}^B=
\left(\begin{matrix}
0&0&0\\
-iU_{\al{\ol{\rho}}}&-V^\beta{}_{\al{\ol{\rho}}}&0\\
-Y_{\ol{\rho}}&-Q_{\ol{\rho}}{}^{\be}&0
\end{matrix}\right),
$$
where the component tensors are given by
$$
S_{\al\beb\rho\sib} = R_{\al\beb\rho\sib} - P_{\al\beb}h_{\rho\sib}
-P_{\rho\sib}h_{\al\beb} - P_{\al\sib}h_{\rho\beb}-P_{\rho\beb}h_{\al\sib}
$$
$$
V_{\al\beb\rho}= \nd_{\beb}A_{\al\rho}+i\nd_{\rho}P_{\al\beb}
-iT_{\rho}h_{\al\beb} -2iT_{\al}h_{\rho\beb}
$$
$$
U_{\al\beb}=\nd_{\al}T_{\beb}+\nd_{\beb}T_{\al}+P_{\al}{}^{\gamma}P_{\gamma\beb}
-A_{\al\gamma}A^{\gamma}{}_{\beb}+Sh_{\al\beb}
$$
$$
Q_{\al\be}=i\nd_0A_{\al\be}-2i\nd_{\be}T_{\al}+2P_{\al}{}^{\gamma}A_{\gamma\be}
$$
$$
Y_{\al}=\nd_0T_{\al}-i\nd_{\al}S+2iP_{\al}{}^{\gamma}T_{\gamma}
-3A_{\al}{}^{\ol{\gamma}}T_{\ol{\gamma}}
$$
$$
V_{\beb\al\sib}={\ol{V_{\beta{\ol{\al}}\sigma}}}, \quad
Q_{{\ol{\al}}\beb}= {\ol{Q_{\al\be}}}, \quad
Y_{\beb}={\ol{Y_{\be}}},
$$
and have the properties:
$$
S_{\al\beb\rho\sib} =
S_{\rho\beb\al\sib}={\ol{S_{\be{\ol{\al}}\sigma{\ol{\rho}}}}},\quad
S_{\al}{}^{\al}{}_{\rho\sib}=0
$$
$$
V_{\al\beb\gamma}=V_{\gamma\beb\al}, \quad V_{\al}{}^{\al}{}_{\gamma}=0
$$
$$
U_{\al\beb}={\ol{U_{\be{\ol{\al}}}}},\quad U_{\al}{}^{\al}=0
$$
$$
Q_{\al\be}=Q_{\be\al}.
$$
This should be compared with the curvature of Chern's Cartan connection
(see \cite{CM}, \cite{web}).
We observe for future reference that the vanishing trace conditions above
imply
\begin{equation}\label{notrace}
\Omega_{\rho}{}^{\rho}{}_A{}^B = 0.
\end{equation}

Next we use the tractor connection to construct a CR invariant second order
differential operator $D$ between tractor bundles.  Let us write
$\ce^*(w,w')$ to indicate any weighted tractor bundle (with arbitrary lists
of conjugated and/or unconjugated upper and/or lower indices).
The {\em tractor D} operator
$$
D_A:\ce^* (w,w')\to \ce_A\otimes \ce^*(w-1,w')
$$
is defined by
$$
D_A f=\left( \begin{array}{c} w(n+w+w') f\\
                            (n+w+w')\nd_\al f\\
  -(\nd^\be\nd_\be f+iw\nd_0 f +w(1+\frac{w'-w}{n+2})Pf)
                            \end{array}\right)
$$
in the realization determined by a choice of $\theta$.  Here
$\nd_\al f$ refers to the tractor connection defined above coupled with
the pseudohermitian connection on the relevant density bundle.  The
$\nd^\be$ in the bottom component additionally couples the pseudohermitian
connection on $\E_\be$.  That this definition is CR invariant is again a
direct calculation of its transformation under change of $\theta$.
The tractor indices on $f$ play no role whatsoever in this calculation
since the tractor connection is invariant.  The first component does not
depend on $\theta$, as required.  For the second component, by
Proposition~\ref{denhat} we have
$\widehat{\nd}_\al f = \nd_\al f + w \Up_\al f$, which is precisely what is
required by the identification.  The third follows upon transforming
the connection using \eqref{vector} and Proposition~\ref{denhat} and
substituting the contraction of \eqref{Ptrans}.  Conjugation produces the
operator
$$
D_{\ol{A}}:\ce^* (w,w')\to \ce_{\ol{A}}\otimes \ce^*(w,w'-1).
$$
Both $D_A$ and $D_{\ol{A}}$ clearly commute with raising and lowering
tractor indices.
A straightforward calculation from the definition shows that one has
\begin{equation}\label{DZ}
D_A Z^A f = (n+w+w'+2)(n+w+1) f
\end{equation}
for $f\in\ce^*(w,w')$,
where here and in what follows we view $Z^A$ as a multiplication operator
and we omit the parentheses when composing operators.

One motivation for the $D_A$ operator is its interpretation in flat
space.  The flat model can be taken to be a real hyperquadric in
$\C{\mathbb P}^{n+1}$.  A section of $\ce(w,w')$ can be interpreted as
a function on the associated affine null cone $\mathcal N$ which is
homogeneous of degree 
$(w,w')$.  If $n+w+w'\neq 0$, such a function $f$ has a homogeneous
extension $\tilde f$
off $\mathcal N$ which is harmonic to first order with respect to the
corresponding Laplacian, and one has
$$
D_A f = (n+w+w')\pa_A {\tilde f}|_{\mathcal N}.
$$
A similar interpretation can be made in the curved case when $M$ is a
hypersurface in $\C^{n+1}$; in this situation the extension is required to
be harmonic to first order with respect to Fefferman's ambient metric.  
Details of this and further connections between tractors and the ambient
construction will be given in \cite{CG2}.

It is often useful to decompose $D_A$ into two pieces.  Given a choice of
$\theta$, let us write $\Box:\ce^*(w,w')\to \ce^*(w-1,w'-1)$ for the
operator
appearing in the bottom slot of the above formula for the tractor $ D$
operator; that is
\begin{equation}\label{boxdef}
\Box f:= \nd^{\al}\nd_{\al} f+iw\nd_0 f +w(1+\frac{(w'-w)}{n+2})Pf,
\end{equation}
and by
$\ol{\Box^{\vphantom{A}}_{\vphantom{A}}}:\ce^*(w,w')\to \ce^*(w-1,w'-1)$
the corresponding conjugate operator given by
$$
\ol{\Box^{\vphantom{A}}_{\vphantom{A}}} f:=
\nd_{\al}\nd^{\al} f-iw'\nd_0 f +w'(1-\frac{(w'-w)}{n+2})Pf.
$$
We also define
$\td_A:\ce^*(w,w')\to\ce_A\otimes\ce^*(w-1,w')$
by
$$
\td_A f:= \left( \begin{array}{c} w f\\
                            \nd_\al f\\
                                  0
 \end{array}\right),
$$
with conjugate
$\td_{\bar{A}}:\ce^*(w,w')\to\ce_{\ol{A}}\otimes\ce^*(w,w'-1)$.
Then clearly we have
\begin{equation}\label{CRexp}
D_A f= (n+w+w')\td_A f-Z_A \Box f.
\end{equation}
A direct calculation shows that as operators on any weighted tractor
bundle,
\begin{equation}\label{boxz}
[\Box,Z_A] =\td_A~.
\end{equation}
Another calculation (using \eqref{notrace}) gives
\begin{equation}\label{boxdiff}
(\Box  - \ol{\Box^{\vphantom{A}}_{\vphantom{A}}}) = (n+w+w')
( i\nd_0+\frac{(w'-w)}{n+2}P)
\end{equation}
on $\ce^*(w,w')$.
The operators $\Box$ and $\td_A$ depend on the choice of $\theta$ and are
not CR invariant.

\section{Invariant Operators Via Tractors} \label{maint}

A first attempt to use tractors to construct an invariant operator from 
densities to densities would be to form $D_AD^A f$.  However, direct
calculation shows that this vanishes identically for $f$ a tractor field of
any weight.
Note, though, that if $f$ has weight $(w,w')$ satisfying
$n+w+w'=0$, then from \nn{CRexp}
$$
D_Af=-Z_A\Box f .
$$
Therefore, $\Box$ is a CR invariant differential operator
$$
\Box : \ce^\ast (w,w')\to \ce^\ast (w-1,w'-1),\qquad n+w+w'=0.
$$
This is a generalization
of the sub-Laplacian of Jerison-Lee (\cite{JLee}), which is
$\Box$ on $\ce(-\frac{n}{2},-\frac{n}{2})$.
The operator
$\ol{\Box^{\vphantom{A}}_{\vphantom{A}}}$ is also CR invariant
when $n+w+w'=0$, but from \eqref{boxdiff} we conclude that
$\Box=\ol{\Box^{\vphantom{A}}_{\vphantom{A}}}$ for such $(w,w')$.

It follows immediately from these observations and the invariance of $D_A$
that for $k\geq 1$, the operator
\begin{equation}\label{tractorform}
D^B\cdots D^A \Box \underbrace{D_A\cdots D_B}_{k-1} \,:\,\,
\ce^*(w,w')\to\ce^*(w-k,w'-k)
\end{equation}
is CR invariant if $n+w+w'=k-1$.  Variants of this operator can also be
formed which are invariant between the same spaces by replacing some of the
indices by barred indices and by reordering
the $D$ factors independently on either side of the $\Box$.  It is not
hard to see that in any choice of CR scale, the principal part of any such
operator agrees with a multiple of that of $\Delta_b^k$, where
$\Delta_b=-(\nd^{\al}\nd_{\al}+\nd^{\ol{\al}}\nd_{\ol{\al}})$.
However, the multiple can vanish.  By calculating these operators
explicitly in flat space, we will determine when the multiple is
nonzero and thereby prove the following strengthening of
Theorem~\ref{intromain}. 
\begin{theorem}\label{main}
For each $(w,w')$ such that $n+w+w'+1=k\in \N$ and
$(w,w') \notin \N_0\times \N_0$,
there is a CR invariant natural differential operator
$$
\P_{w,w'}: \ce^*(w,w')\to \ce^*(w-k,w'-k)
$$
whose principal part agrees with that of $\Delta_b^k$.
\end{theorem}
\noindent
For fixed $(w,w')$
as in Theorem~\ref{main} we shall in general obtain 
several operators $\P_{w,w'}$ corresponding to different choices of barred
and unbarred indices and different orderings of the $D$ factors.  We shall
see that these operators all agree for a flat CR structure.

It is straightforward but tedious to check directly from the definition
that for a flat pseudohermitian structure,
\begin{equation}\label{commute}
[D_B,D_C]=0, \quad [D_B,D_{\bar{C}}]=0,
~\mbox{ and }~ [D_{\bar{B}},D_{\bar{C}}]=0
\end{equation}
as operators on any weighted tractor bundle.  Since $D_A$ is
CR invariant, this remains true if the structure is only CR flat.  So in
the flat case, it is clear that the order of the $D$ factors is irrelevant.

The following proposition calculates the operators in flat space.  The
proof of Theorem~\ref{main} only uses the cases $k_1=0$ or $k_2=0$, but we
include the general case for completeness.

\begin{proposition}\label{flatgoody}  Let $(w,w')$ satisfy
$n+w+w'+1 =k \in \N$. Let $k_1, k_2 \in \N_0$ be such that
$k_1+k_2=k-1$.  For a flat pseudohermitian
structure, we have as operators on $\ce^*(w,w')$:
\begin{equation}\label{prelim}
\Box\underbrace{D_B\cdots D_I}_{k_1}\underbrace{D_{\bar{J}}\cdots
  D_{\bar{Q}}}_{k_2}
=(-1)^{k-1}
Z_B\cdots Z_I Z_{\bar{J}}\cdots Z_{\bar{Q}} \,\Box^k
\end{equation}
\begin{equation}\label{flatop}
\begin{split}
D^{\bar{Q}}\cdots D^{\bar{J}}D^I\cdots  D^B & \Box
\underbrace{D_B\cdots D_I}_{k_1}
\underbrace{D_{\bar{J}}\cdots D_{\bar{Q}}}_{k_2} \\
& = \,(-1)^{k-1}(k-1)!\prod_{i=0}^{k_1-1} (w-i)\prod_{j=0}^{k_2-1}
(w'-j)\,\Box^k
\end{split}
\end{equation}
\begin{equation}\label{conjugation}
%\Box^k = \ol{\Box^{\vphantom{A}}_{\vphantom{A}}}{}^k
\Box^k = \ol{\Box}^k.
\end{equation}
\end{proposition}
\medskip
\noindent
We have written $\Box^k$ for the
$k$-fold composition of $\Box$ with itself.  Note that each factor is
acting on a different density space, so each $\Box$ in the composition is
really a different operator.  In case $k_1$ or $k_2$ equals $0$, the empty
product in \eqref{flatop} is to be interpreted as $1$.

\medskip
\noindent
{\em Proof of Theorem~\ref{main} using Proposition~\ref{flatgoody}.}
Conjugating if necessary, we may assume without loss of generality that
$w\notin \N_0$.
Set $k_1=k-1$, $k_2=0$.  The numerical factor on the right hand side of
\eqref{flatop} is then nonzero.  Since the principal
part of $-2\Box$ agrees with that of $\Delta_b$, 
we may define $\P_{w,w'}$ by multiplying the left hand side 
of \eqref{tractorform} by the appropriate constant.  It is clear that 
$\P_{w,w'}$ is a CR invariant operator, and 
expanding \eqref{tractorform} shows that it is natural in the sense that  
it is given by a universal formula in terms of the pseudohermitian
metric, curvature, torsion, and connection.  In 
flat space it has the correct principal part, and it is easily seen
that the correction terms in curved space are of lower order.
\stopthm

We remark that the proof shows that in Theorem~\ref{main} a stronger
statement can be made
about the principal part of $\P_{w,w'}$:  its nonisotropic principal
part, in which derivatives with respect to $T$ are weighted by a factor of 
2, agrees with that of $(-2\Box)^k$.

The proof of Proposition~\ref{flatgoody} is preceeded by three lemmas.

\begin{lemma}\label{one}
Let $k_1, k_2 \in \N_0$ and let
$(w,w')$ satisfy $n+w+w'=k_1+k_2$.  For a CR flat 
structure, we have as operators on $\ce^*(w,w')$:
\[
\Box\underbrace{D_B\cdots D_I}_{k_1}\underbrace{D_{\bar{J}}\cdots
  D_{\bar{Q}}}_{k_2}
=(-1)^{k_1+k_2}
Z_B\cdots Z_I Z_{\bar{J}}\cdots Z_{\bar{Q}} \,\bbox_{k_1,k_2}
\]
for a differential operator $ \bbox_{k_1,k_2} $.
\end{lemma}
\begin{proof}
It follows from \eqref{CRexp} that if $f\in \ce^*(w,w')$, then
$$
Z_{A} \Box D_{B}D_C \cdots D_I  D_{\bar{J}}\cdots
  D_{\bar{Q}} f = -D_A D_{B}D_C \cdots D_I  D_{\bar{J}}\cdots
  D_{\bar{Q}} f.
$$
Skewing on $A$ and $B$ and recalling \eqref{commute} gives
$$
Z_{[A} \Box D_{B]}D_C \cdots D_I  D_{\bar{J}}\cdots
  D_{\bar{Q}} f = 0.
$$
Again using \eqref{commute}, we can commute the $D$'s to conclude
that the commutator can be taken on
any of the unbarred indices:
$$
Z_{[A}\Box D_{|B}\cdots D_{E|}D_{F]}D_G\cdots D_I  D_{\bar{J}}\cdots
  D_{\bar{Q}} f=0
$$
for $F=B,C,\cdots I$.  Recalling \eqref{boxdiff}, we obtain similarly
$$
Z_{[\bar{A}}\Box D_{|B}\cdots D_I D_{\bar{J}}\cdots
D_{\bar{M}|}D_{\bar{N}]}D_{\bar{P}}\cdots
D_{\bar{Q}} f=0
$$
for $\bar{N}=\bar{J}\cdots \bar{Q}$.
Therefore
$$
\Box D_{B}D_C \cdots D_I  D_{\bar{J}}\cdots
  D_{\bar{Q}} f
$$
must be proportional to each of $Z_B, \cdots, Z_{\bar{Q}}$, and the result
follows.
\end{proof}

\begin{lemma}\label{two}
In Lemma~\ref{one}, we have
$\bbox_{k_1,k_2} = \bbox_{k_1',k_2'} $ if
$k_1+k_2=k_1'+k_2'$.
\end{lemma}
\begin{proof}
It suffices to show that for $k_2\geq 1$ we have
$\bbox_{k_1,k_2} = \bbox_{k_1+1,k_2-1}$ as operators on
$\ce^*(w,w')$ with $(w,w')$ as in Lemma~\ref{one}.
This is a consequence of \eqref{CRexp}, \eqref{commute} and \eqref{boxdiff}
as follows:
\[
\begin{split}
(-1)^{k_1+k_2}Z_AZ_B \cdots Z_I &  Z_{\bar{J}}\cdots Z_{\bar{Q}}\,
\bbox_{k_1,k_2} \\
&=Z_A\Box\underbrace{D_B\cdots D_I}_{k_1}\underbrace{D_{\bar{J}}\cdots
  D_{\bar{Q}}}_{k_2} \\
&=- D_A D_B\cdots D_I D_{\bar{J}}D_{\bar{K}}\cdots
  D_{\bar{Q}} \\
&=- D_{\bar{J}}D_A D_B\cdots D_I D_{\bar{K}}\cdots
  D_{\bar{Q}} \\
&= Z_{\bar{J}}\Box D_A D_B\cdots D_I D_{\bar{K}}\cdots
  D_{\bar{Q}}  \\
&=(-1)^{k_1+k_2} Z_AZ_B\cdots Z_I Z_{\bar{J}}\cdots Z_{\bar{Q}}\,
\bbox_{k_1+1,k_2-1}.
\end{split}
\]
\end{proof}

\begin{lemma}\label{three}
For flat pseudohermitian structures, we have
$$
[\Box,\td_A]=0, \qquad [\Box^k,Z_A] = k\Box^{k-1}\td_A
$$
on any weighted tractor bundle.
\end{lemma}
\begin{proof}
The first equation follows by direct calculation from the definitions.
The second is a consequence of the first together with \eqref{boxz}.
\end{proof}
\medskip
\noindent
{\em Proof of Proposition~\ref{flatgoody}.}
According to Lemmas \ref{one} and \ref{two}, in order to establish
\eqref{prelim} it suffices to show that for $k\geq 1$, we have for flat
pseudohermitian structures
$\bbox_{k-1,0}=\Box^k$ as operators on $\ce^*(w,w')$ for $n+w+w'+1=k$.
We prove this by induction on $k$.  The case $k=1$ is clear.
Suppose the statement is true for $k$ and let $f\in \ce^*(w,w')$ with
$n+w+w'=k$.  Note that it follows from Lemma~\ref{one} that
$$
Z_A \bbox_{k,0} f =-\bbox_{k-1,0} D_A f.
$$
Combining this with the induction hypothesis and then using \eqref{CRexp}
and Lemma~\ref{three} gives
$$
Z_A \bbox_{k,0} f =-\Box^k D_A f = -\Box^k(k\td_A -Z_A\Box)f
= Z_A \Box^{k+1}f,
$$
as desired.

Now \eqref{flatop} follows upon repeatedly applying \eqref{DZ} and its
conjugate to \eqref{prelim}, and \eqref{conjugation} follows upon comparing
\eqref{prelim} for $k_2=0$ with the conjugate of \eqref{prelim} for
$k_1=0$, recalling from \eqref{boxdiff} that $\Box = {\ol{\Box}}$ on
$\ce^*(w,w')$ for $n+w+w'=0$.
\stopthm

For CR flat structures, Lemmas \ref{one} and \ref{two} show that the
operator $\bbox_{k_1,k_2}$ is CR invariant and depends only on 
$k_1+k_2$, and its principal part is easily determined.  It follows
that for CR flat structures we can extend the definition of $\P_{w,w'}$ to
the case $n+w+w'+1=k\in \N$ without the restriction $k\leq n+1$ 
by setting $\P_{w,w'}=(-2)^k\bbox_{k_1,k_2}$.  Then for CR flat structures 
$\P_{w,w'}$ is 
independent of the choice of barred and unbarred indices (and
of their ordering, as already observed).

In the curved case,
one expects by analogy with the conformal case (see the discussion in
\cite{Gosrni}, \cite{Esrni}) that it is no longer generally true that
$\Box D_B\cdots D_I D_{\bar{J}}\cdots
  D_{\bar{Q}}$ is proportional to $Z_B\cdots Z_I Z_{\bar{J}}\cdots
Z_{\bar{Q}}$, necessitating
the application of $D^{\bar{Q}}\cdots D^{\bar{J}}D^I\cdots  D^B$ as in
\eqref{flatop} and leading to the
condition $(w,w') \notin \N_0\times \N_0$ in Theorem~\ref{main}.
However, for $k=2$ this is
not necessary even in the curved case:  if $f\in \ce(w,w')$ with
$n+w+w'=1$, one
can see by direct calculation from the definitions that the first two slots
of $\Box D_A f$ vanish, so we have 
$4 \Box D_A f = -Z_A L_{w,w'}f$ for a fourth 
order operator $L_{w,w'}$ whose principal part agrees with
that of $\Delta_b^2$.  Applying $D^A$ shows that $L_{w,w'}= \P_{w,w'}$ if 
$w\neq 0$, so in this case nothing new is obtained.
However, if $n=1$ we conclude the existence of an 
invariant operator $L_{0,0}:\ce(0,0)\to\ce(-2,-2)$ not covered by
Theorem~\ref{main}.  For $n=1$, we may therefore extend the definition of
the family $\P_{w,w'}$ by setting $\P_{0,0}=L_{0,0}$.  We shall
recover $\P_{0,0}$ and show the 
existence of higher
dimensional analogues mapping $\ce(0,0)\to \ce(-n-1,-n-1)$ in the next
section using the ambient metric.

Next we show that as long as the $D$'s are ordered consistently, the
operators $\P_{w,w'}$
constructed in Theorem~\ref{main} are self-adjoint.  Observe first that
under the change 
of scale $\hat{\theta}=e^{\Up}\theta$, the volume form
$\theta\wedge(d\theta)^n $ multiplies by $e^{(n+1)\Up}$, so the bundle of
volume densities can be canonically identified with $\ce(-n-1,-n-1)$ and
$\int_M f$ is invariantly defined for $f \in \ce(-n-1,-n-1)$.
(Throughout this discussion we assume that all sections are compactly
supported.)  In general,
if $L:E\to F$ is a differential operator between complex vector bundles $E$
and $F$, then the formal adjoint $L^*:{\ol{F}}^*\otimes \ce(-n-1,-n-1)\to
{\ol{E}}^*\otimes \ce(-n-1,-n-1)$ is
defined by
$$
\int_M \langle Lu,{\ol{v}}\rangle = \int_M \langle u ,{\ol{L^*v}}\rangle 
$$  
for $u \in \Gamma(E)$, $v\in \Gamma({\ol{F}}^*\otimes \ce(-n-1,-n-1))$, 
where ${\ol{E}}^*$ denotes the conjugate of the dual bundle to $E$.
We have ${\ol{\ce(w,w')}}^* = \ce(-{\ol{w'}},-{\ol{w}})$.  If
$n+w+w'+1=k\in \N$, then by our implicit assumption that $w-w'\in \mathbb
Z$ it follows that $w,w'\in \mathbb R$, so the operators
$\P_{w,w'}:\ce(w,w')\to \ce(w-k,w'-k)$ constructed in Theorem~\ref{main}
have the property that also
$\P_{w,w'}{}^*:\ce(w,w')\to \ce(w-k,w'-k)$.  Recall that the operators
$\P_{w,w'}$ depend on a choice of ordering of barred and/or unbarred
indices before and after the middle $\Box$.  We shall say that the indices
are ordered consistently if the order of the indices before the $\Box$
is the opposite of those after the $\Box$.
\begin{proposition}\label{duality}
If $(w,w')$ is as in Theorem~\ref{main} and if the indices are ordered
consistently, then 
$\P_{w,w'}:\ce(w,w')\to \ce(w-k,w'-k)$  is self-adjoint.
\end{proposition}
\begin{proof}
The proposition clearly follows if we show that $\Box$ is
self-adjoint on $\ce^*(w,w')$ for $n+w+w'=0$ and that
\begin{equation}\label{Diby}
\int_M D_A f \cdot g^A = \int_M f\cdot D_Ag^A
\end{equation}
for $f\in \ce^*(w,w')$, $g^A\in \ce^A\otimes \ce^*(-n-w,-n-w'-1)$, where
the supressed tractor indices on $f$ and $g^A$ are the same and we
have denoted by `$\cdot$' the full contraction of these tractor
indices.

An easy calculation shows that the
divergence operator $f_{\al}\to \nd^{\al}f_{\al}$ is CR invariant
$:\ce_{\al}(-n,-n)\to \ce(-n-1,-n-1)$ and
integration by parts (see (2.18) of \cite{Lee2}) shows that
\begin{equation}\label{firstiby}
\int_M \nd^{\al}f_{\al} = 0
\end{equation}
for $f_\al \in \ce_{\al}(-n,-n)$. 
Also it is easily seen from \eqref{boxdef}, \eqref{dencomm} and
\eqref{notrace} that
$$
n\Box  = n\nd^{\al}\nd_{\al}+w[\nd^{\al},\nd_{\al}] +w(n+w-w')P
$$
on $\ce^*(w,w')$.  Therefore for $f\in\ce^*(w,w')$, $g\in\ce^*(-n-w,-n-w')$
we have
\begin{equation}\label{boxiby}
\begin{split}
n\int_M \Box f &\cdot g \\
&=\int_M (n\nd^{\al}\nd_{\al}f+w[\nd^{\al},\nd_{\al}]f +w(n+w-w')Pf)\cdot
g\\
&=\int_M f\cdot(n\nd_{\al}\nd^{\al}g +w[\nd_{\al},\nd^{\al}]g
+w(n+w-w')Pg)\\
&=\int_M f\cdot(n\nd^{\al}\nd_{\al}g -(n+w)[\nd^{\al},\nd_{\al}]g
+w(n+w-w')Pg)\\
&=n\int_M f\cdot (\Box g +(n+w+w')Pg).
\end{split}
\end{equation}
Self-adjointness of $\Box$ when $n+w+w'=0$ follows immediately from this
together with the fact that $\Box = {\ol{\Box}}$ when $n+w+w'=0$.
Finally, if $f\in \ce^*(w,w')$ and
$$
g_A=\left( \begin{array}{c} \sigma\\
                            \tau_{\al}\\
                             \rho \end{array}\right)
\in \ce_A\otimes \ce^*(-n-w,-n-w'-1),
$$
then the definition of $D_A$ gives
$$
D_Af\cdot g^A = w(n+w+w')f\cdot \rho +(n+w+w')\nd_{\al}f\cdot\tau^{\al}
-\Box f \cdot \sigma,
$$
while a computation shows that
$$
D_Ag^A=w(n+w+w')\rho -(n+w+w')\nd_{\al}\tau^{\al}
-(\Box + (n+w+w')P)\sigma.
$$
Therefore \eqref{Diby} follows from \eqref{firstiby} and \eqref{boxiby}.
\end{proof}

\section{Invariant Operators Via Fefferman Metric} \label{amb}
In this section we review Lee's formulation (\cite{Lee1}) of
the Fefferman conformal structure and the
construction in \cite{GJMS} of the invariant powers of the Laplacian
on a conformal manifold via the ambient metric.  Combining these produces
CR invariant powers of the sub-Laplacian.  For a hypersurface in
$\C^{n+1}$, this construction can be expressed in 
terms of Fefferman's original formulation, and reexpressed in terms of the
associated Cheng-Yau metric.     

In Lee's formulation (see also \cite{Fa}, \cite{BDS}), the Fefferman metric
of a CR manifold $M$ lives on the circle bundle $\K^*/\R_+$, where 
$\K$ denotes the canonical bundle of $M$, and $\K^*=\K\setminus 0$.  A  
section $\zeta$ of $\K^*$ determines a fiber variable $\psi$ on 
$\K^*/\R_+$ by the requirement that $e^{i\psi}\zeta$ be in the given
$\R_+$-equivalence class.
Suppose that $\theta$ is a choice of
pseudohermitian form, that $\zeta$ is volume normalized with respect to
$\theta$, and choose an admissible coframe $\theta^{\alpha}$ such that 
$\zeta = \theta\wedge\theta^1\wedge\cdots\wedge\theta^n$.  The 1-form
$\sigma$ on $\K^*/\R_+$ defined by 
$$
(n+2) \sigma = d\psi + i\omega_{\al}{}^{\al} - \frac{1}{2(n+1)}R\theta 
$$
is independent of the choice of $\zeta$ and $\theta^{\al}$, so depends only
on $\theta$ and is globally
defined.  The Fefferman metric associated to $\theta$ is the
metric of signature $(2p+1,2q+1)$ given by
$$
g=h_{\al\beb}\theta^{\al}\cdot\theta^{\beb} +2\theta\cdot\sigma.
$$
In \cite{Lee1} it is shown that if $\hat{\theta}=e^{\Up}\theta$, then
$\hat{g}=e^{\Up}g$, so that the conformal class of $g$ is CR invariant. 
Lee also explicitly calculated the connection forms and Ricci curvature of
$g$.  From his expression for the connection, it follows that all
components of the curvature tensor of $g$ and its iterated covariant
derivatives are given by universal expressions in the pseudohermitian
metric, curvature and torsion and their pseudohermitian covariant
derivatives.  

As in Section 2, we assume that $\K$ admits an 
$(n+2)^{\rm nd}$ root, which we fix.  For our purposes it will be useful to
work on the circle bundle  
$\Co=(\K^*)^{1/(n+2)}/\R_+$.  The Fefferman metric pulls back via the 
$(n+2)^{\rm nd}$ power map to a metric on $\Co$
which we shall also denote by $g$.  A fiber variable $\g$ on $\Co$
satisfies $(n+2)\g = \psi$.

The metric bundle of a manifold $\Co$ of dimension $N$ with a conformal
class of metrics $[g]$ of signature $(P,Q)$ is the ray subbundle 
${\mathcal G}\subset S^2T^*\Co$ of multiples of the metric:  if $g$ is a
representative metric, then the fiber of ${\mathcal G}$ over $p\in \Co$ is
$\{t^2g(p):t>0\}$. The bundle of conformal densities of weight $w\in \C$ is 
$\ce ( w ) = {\mathcal G}^{-w/2}$,
where by abuse of notation we have denoted by ${\mathcal G}$ also the line
bundle associated to the ray bundle defined above.  The main result of
\cite{GJMS} is the existence, for $k\in \N$ satisfying $k\leq N/2$ if $N$
is even, of a conformally invariant natural differential operator 
$P_k:\ce(-N/2+k)\rightarrow \ce(-N/2-k)$ with principal part equal to that
of $\Delta^k$.  
These operators are constructed in \cite{GJMS} using the ambient metric
of \cite{FGast}.  Denote by 
$\pi:\mathcal G \rightarrow \Co$ the natural projection of the metric
bundle, and by ${\bf g}$ the tautological  
symmetric 2-tensor on $\mathcal G$ defined for $(p,g)\in \mathcal G$
and $X,Y\in T_{(p,g)}{\mathcal G}$ 
by ${\bf g}(X,Y)=g(\pi_*X,\pi_*Y)$.  There are dilations 
$\delta_s:{\mathcal G}\rightarrow {\mathcal G}$ for $s>0$ given by
$\delta_s(p,g)=(p,s^2g)$, and we have  
$\delta_s^*{\bf g} = s^2{\bf g}$.  Denote by $S$ the infinitesimal dilation
vector field $S=\frac{d}{ds} \delta_s |_{s=1}$.  Define the ambient space 
${\tilde {\mathcal G}}={\mathcal G}\times (-1,1)$.  Identify ${\mathcal G}$
with its image under 
the inclusion $\iota:{\mathcal G}\rightarrow {\tilde {\mathcal G}}$ given
by $\iota(g)=(g,0)$ for $g\in {\mathcal G}$.  The dilations $\delta_s$ and
infinitesimal generator $S$ extend naturally to ${\tilde {\mathcal G}}$.  
The ambient metric ${\tilde g}$ is a metric of signature $(P+1,Q+1)$ on 
${\tilde {\mathcal G}}$
which satisfies the initial condition
$\iota^* {\tilde g}={\bf g}$, is homogeneous in the sense that 
$\delta_s^* {\tilde g}=s^2{\tilde g}$, and is an asymptotic solution of 
$\Ric({\tilde g})=0$ along $\mathcal G$.  For $N$ odd, these conditions
uniquely 
determine a formal power series expansion for ${\tilde g}$ up to
diffeomorphism, but for $N$ even and $N>2$, a formal power series solution
exists in general only to order $N/2$.   

An element of $\ce(w)$ can be regarded as a homogeneous function of degree
$w$ on $\mathcal G$.  It is shown in \cite{GJMS} that the same operator
$P_k$ arises in two ways:
\begin{enumerate}
\item By extending a density $f\in \ce(-N/2+k)$ to a function
${\tilde f}$ homogeneous of degree $-N/2+k$ on ${\tilde {\mathcal
G}}$, applying ${\tilde \Delta}^k $, where $\tilde \Delta$ denotes the
Laplacian in the metric ${\tilde g}$, and restricting back
to $\mathcal G$:  $P_k f = {\tilde \Delta}^k \tilde f |_{\mathcal G}$.
\item As the normalized obstruction to extending 
${\tilde f}\in \ce(-N/2+k)$ to a smooth function ${\tilde F}$ 
homogeneous of degree $-N/2+k$ on ${\tilde
{\mathcal G}}$, such that ${\tilde F}$ satisfies
${\tilde \Delta} {\tilde F}= 0$ to infinite order along $\mathcal G$. 
\end{enumerate}
For $N$ even, the condition $k\leq N/2$ is needed to ensure that the
construction does not involve too many derivatives of the ambient metric. 
We remark that 2.\ above can be restated in terms of non-smooth solutions:
there are infinite order formal solutions whose expansion includes a 
$\log$ term, and the invariant
operator applied to $f$ can be characterized as a multiple of the
restriction to $\mathcal G$ of the coefficient of the $\log$ term of a
formal solution which equals $f$ on $\mathcal G$.

Let now $M$ be a CR manifold of dimension $2n+1$.
We associate to $M$ its Fefferman conformal
manifold $(\Co,[g])$, of dimension $N=2n+2$.  Volume normalization gives a
canonical identification between $(\K^*)^{1/(n+2)}$ (with usual scalar
multiplication) and the metric bundle of 
$(\Co,[g])$ (with dilations $\delta_s$) as 
$\R_+$-principal bundles over $\Co$.  Now
a CR density $f \in \ce(w,w')$ 
may be viewed as a smooth function on $(\K^*)^{1/(n+2)}$ homogeneous of
degree $(w,w')$ in the sense that 
$f(\lambda \xi) = \lambda^w{\overline{\lambda}}{}^{w'}f(\xi)$ for 
$\xi\in (\K^*)^{1/(n+2)}$.  
Therefore $\ce(w,w')$ may be regarded as the subspace of the space of
conformal densities $\ce(w+w')$ which satisfy 
\begin{equation}\label{densityinclusion}
(e^{i\phi})^*f = e^{i(w-w')\phi}f\,, \qquad \phi \in \R,
\end{equation}
where on the left hand side,
$(e^{i\phi})^*$ denotes pull back by the isometry of $(\Co,g)$ given by 
multiplication by $e^{i\phi}$.

\medskip
\noindent
{\em Proof of Theorem~\ref{ambient}.}
The conformally invariant operator $P_k$ of \cite{GJMS} satisfies 
$P_k:\ce(w+w')\rightarrow \ce(w+w'-2k)$.  Since multiplication by
$e^{i\phi}$ is an isometry of $(\Co,g)$ and $P_k$ is a natural operator, it 
follows that $P_k f$ satisfies \eqref{densityinclusion} if $f$ does.  
Therefore $P_k$ induces an operator from $\ce(w,w')$ to $\ce(w-k,w'-k)$. 
Define $P_{w,w'} = 2^{-k}P_k|_{\ce(w,w')}$.  The CR invariance of
$P_{w,w'}$ follows from the conformal invariance of $P_k$.  Using the
formulae in \cite{Lee1} for the connection of the Fefferman metric, 
one can see
without difficulty that any component of an
iterated covariant derivative of $f\in \ce(w,w')$ with respect to the
Levi-Civita connection of $g$ has a 
universal expression in terms of pseudohermitian covariant derivatives of
$f$ and the pseudohermitian metric, torsion and curvature.  Since $P_k$ is
a natural 
differential operator associated to pseudo-Riemannian manifolds, it follows 
that $P_{w,w'}$ is a natural differential operator associated to
pseudohermitian manifolds in the sense that it also has such a universal 
expression.  It is
easily seen that the principal part of $P_{w,w'}$ agrees with that of
$\Delta_b^k$. 
\stopthm

For $k=1$ it is straightforward to carry out the calculation 
of $P_{w,w'}$ from the conformal Laplacian $P_1$; one finds in that case
that $P_{w,w'}=\P_{w,w'}=-2\Box$.

The conformal operators $P_k$ have been shown in \cite{GZ} and \cite{FG2}
to be self-adjoint; the argument in \cite{FG2} is formal and valid in
general signature.  From this, the self-adjointness of $P_{w,w'}$ is an
easy consequence:

\begin{proposition}\label{ambientsa}
If $n+w+w'+1=k\leq n+1$, then the operator $P_{w,w'}:\ce(w,w')\rightarrow
\ce(w-k,w'-k)$ is self-adjoint.
\end{proposition}
\begin{proof}
If $f_1, f_2 \in \ce(w,w')$, then $P_{w,w'}f_1\,\,{\ol{f_2}} 
=2^{-k}  P_kf_1\,\,{\ol{f_2}}  \in \ce(-n-1,-n-1) \subset \ce(-N)$ is invariant
under rotation in $\Co$.  It follows that $\int_M  P_{w,w'}f_1 \,\,{\ol{f_2}}$
is a constant multiple of $\int_{\Co}  P_kf_1\,\,{\ol{f_2}}$, so the
self-adjointness of $P_{w,w'}$ follows from that of $P_k$.
\end{proof}

A nondegenerate hypersurface $M \subset \C^{n+1}$ inheirits a CR
structure.  This was the setting for Fefferman's original construction of
the ambient metric and conformal metric in \cite{f}, \cite{F2}.  
The ambient metric is defined above on $\mathcal G_{\Co}\times (-1,1) = 
(\K_M^*)^{1/(n+2)} \times (-1,1)$, which we may identify with 
$(\K_{\C^{n+1}}^*)^{1/(n+2)}$ over a neighborhood $U\subset \C^{n+1}$ of
$M$, since $\K_M$ is the restriction of $\K_{\C^{n+1}}$.  
Introduce a variable $z^0\in \C$, which is interpreted
as a fiber coordinate on $(\K_{\C^{n+1}})^{1/(n+2)}$ relative to the
trivialization defined by
$(dz^1\wedge\ldots\wedge dz^{n+1})^{1/(n+2)}$.  So the ambient metric is
defined on $\C^*\times U$.  Fefferman
showed that $M$ has a smooth defining function $\phi$ which
solves $J(\phi)= 1 + O(\phi^{n+2})$, where 
$$
J(\phi)=(-1)^{p+1} \det \left(\begin{matrix}
\phi&\phi_{\bb} \\
\phi_{a}&\phi_{a\bb}
\end{matrix}\right)_{1\leq a,b \leq n+1}
$$
Here the subscripts denote coordinate derivatives and the Levi form
$-\phi_{a\bb}|_{T^{1,0}M}$ has signature $(p,q)$, $p+q=n$.  
The function
$Q:=-|z^0|^2\phi(z)$ satisfies
$
(-1)^{q+1}\det ( \pa^2_{A\Bb}Q)=|z^0|^{2(n+1)}J(\phi),
$
where $A,B=0,\ldots,n+1$.  The Hermitian matrix
$(\pa^2_{A\Bb}Q)$ is therefore nondegenerate and
defines a K\"ahler metric ${\tilde g}$ of signature $(p+1,q+1)$
which is approximately Ricci 
flat along $M$.  This metric is clearly homogeneous 
of degree 2 and it is shown in \cite{Lee1} that its restriction to $S^1
\times M$ (Fefferman's original definition of the Fefferman metric) agrees
with the definition above for the Fefferman  
metric of the pseudohermitian form $\theta = i(\pa-{\ol{\pa}})\phi/2$. 
Therefore the K\"ahler metric ${\tilde g}$ is the ambient metric associated
to $M$.   

We remark that in this formulation, the pseudohermitian form 
$\theta = i(\pa-{\ol{\pa}})\phi/2$ is not general:  the condition 
$J(\phi)=1$ on $M$ means that the closed form $dz^1\wedge\ldots\wedge
dz^{n+1}$ is volume normalized with respect to $\theta$.  This 
forces conditions on $\theta$:  if $n>1$ then $\theta$ is pseudo-Einstein 
(\cite{Lee2}); the corresponding condition for $n=1$ is given in
\cite{hirachi}.  Invariance under rescaling $\theta$ is
replaced by invariance under biholomorphic changes of coordinates.

For $M$ a hypersurface in $\C^{n+1}$, a section of $\ce(w,w')$ can be
regarded as a function on $\C^*\times M$ homogeneous of degree $(w,w')$ in
$z^0$.  The invariant operators $P_{w,w'}$ are therefore
realized as 
the obstruction to extending such a function to be homogeneous and harmonic
with respect to the Laplacian ${\tilde{\Delta}}$ in the 
K\"ahler metric ${\tilde g}$, or equivalently by applying
${\tilde{\Delta}}^k$ to an arbitrary homogeneous extension.  

We next show how this realization of the operators $P_{w,w'}$ can be
reformulated in terms of the K\"ahler metric
\begin{equation}\label{poincare}
g_{a\bb}=(\log\phi^{-1})_{a\bb}
\end{equation}
on $\{\phi > 0\}\subset\C^{n+1}$.  If $\phi$ solves $J(\phi)=1$ exactly,
then $g$ is the K\"ahler-Einstein metric of Cheng-Yau \cite{chengyau}. 
We begin by considering the relationship between the Laplacians of $\gt$
and $g$ for a general metric of the form \eqref{poincare}, not necessarily
assuming that $\phi$ is an approximate solution of $J(\phi)=1$.
We continue to use capital indices $A,B$ which run between $0$ and $n+1$,
so for example
\begin{equation}\label{ambmetric}
\gt_{A\Bb} =- \left(
\begin{matrix}
\phi&\zb^0\phi_{\bb}\\
z^0\phi_a&|z^0|^2\phi_{a\bb}\\
\end{matrix}
\right).
\end{equation}
Lower case indices $a,b$ will be raised and lowered using
\begin{equation}\label{ke}
g_{a\bb}=-\phi^{-1}\phi_{a\bb}+\phi^{-2}\phi_a\phi_{\bb}
\end{equation}
and its inverse.

As observed in \cite{leemel}, the choice of
defining function $\phi$ determines near $M$ a $(1,0)$ vector field $\xi^a$
and a scalar function $r$ (called the transverse curvature in
\cite{GL}), as follows.  Observe first that since $\phi_{a\bb}|_{T^{1,0}M}$
is nondegenerate, there is a unique direction in $T^{1,0}\C^{n+1}$
which is orthogonal to $\ker \pa\phi$ with respect to $\phi_{a\bb}$.
Therefore, near $M$ there is a unique $(1,0)$ vector field $\xi^a$
satisfying
\begin{equation}\label{def}
\phi_{a\bb}\xi^{\bb} = r \phi_a, \qquad  \xi^a\phi_a = 1
\end{equation}
for some uniquely determined smooth function $r$, and we have
\begin{equation}\label{r}
r=\phi_{a\bb} \xi^a\xi^{\bb}.
\end{equation}
\begin{proposition}\label{laplprop}
When acting on functions which are homogeneous of degree $(w,w')$ with
respect to $z^0$, we have
\begin{equation}\label{homoglapl}
(\phi|z^0|^2)\Dt=\Delta_g + \frac{w'\phi}{1-r\phi}\xi^a\pa_a
+ \frac{w\phi}{1-r\phi}\xi^{\bb}\pa_{\bb} - \frac{ww'r\phi}{1-r\phi}.
\end{equation}
Here ${\tilde \Delta}$ denotes the K\"ahler Laplacian 
$-{\tilde g}^{A\Bb}\pa^2_{A\Bb}$ and similarly for $\Delta_g$.
\end{proposition}
\begin{proof}
Our first observation is a matrix factorization of $\gt$ given by
\eqref{ambmetric}:
\begin{equation}\label{factor}
\gt_{A\Bb}=
\left(
\begin{matrix}
\phi^{1/2}&0\\
z^0\phi^{-1/2}\phi_a&z^0\phi^{1/2}\delta_a{}^c
\end{matrix}
\right)
\left(
\begin{matrix}
-1&0\\
0&g_{c\ol{d}}
\end{matrix}
\right)
\left(
\begin{matrix}
\phi^{1/2}&\zb^0\phi^{-1/2}\phi_{\bb}\\
0&\zb^0\phi^{1/2}\delta^{\ol{d}}{}_{\bb}
\end{matrix}
\right).
\end{equation}
Taking inverses and multiplying out gives
\begin{equation}\label{inverse}
(|z^0|^2\phi) \gt^{A\Bb}=
\left(
\begin{matrix}
|z^0|^2(\phi^{-2}\phi^a\phi_a - 1)&-z^0\phi^{-1}\phi^{\bb}\\
-\zb^0\phi^{-1}\phi^a&g^{a\bb}
\end{matrix}
\right),
\end{equation}
from which we obtain
\begin{equation}\label{laplacian}
(|z^0|^2\phi) \Dt
 = \Delta_g
+\phi^{-1}\phi^a\pa_a\zb^0\pa_{\ol{0}}
+\phi^{-1}\phi^{\bb}\pa_{\bb}z^0\pa_0
-(\phi^{-2}\phi^a\phi_a-1)|z^0|^2\pa^2_{0\ol{0}}.
\end{equation}
We want to rewrite the coefficients of the last three terms.  Define 
$l_{a\bb}=r\phi_a\phi_{\bb}-\phi_{a\bb}$ so that $l_{a\bb}\xi^{\bb}=0$.
Then $l_{a\bb}$ is nondegenerate on $\ker \pa \phi$, so we may consider its
inverse there, extend to annhilate $\phi_a$, and thus obtain a tensor
$k^{a\bb}$ defined by the relations $k^{a\bb}\phi_{\bb}=0$ and
$k^{a\bb}l_{c\bb}= \delta^a{}_c - \xi^a\phi_c$.  Now \eqref{ke} gives
$g_{a\bb}=\phi^{-1}(l_{a\bb}+(1-r\phi)\phi^{-1}\phi_a\phi_{\bb})$, so
$g^{a\bb} = \phi(k^{a\bb}+\phi(1-r\phi)^{-1}\xi^a\xi^{\bb})$.  {From} this
we conclude that $\phi^{-1}g^{a\bb}$ extends smoothly across $\phi=0$, and
that
\begin{equation}\label{raise}
\phi^a = \frac{\phi^2}{1-r\phi}\xi^a,\qquad \phi^a\phi_a =
\frac{\phi^2}{1-r\phi}.
\end{equation}
Substituting into \eqref{laplacian} and using homogeneity gives
\eqref{homoglapl}.
\end{proof}
\noindent
{\bf Remarks.\,\,}
The factorization \eqref{factor} also gives a direct relation between the
metrics $\gt$ and $g$:
$\,\,
\gt = (|z^0|^2\phi)\left[ g-(\frac{dz^0}{z^0}+\frac{\pa\phi}{\phi})
(\frac{d\zb^0}{\zb^0}+\frac{\ol{\pa}\phi}{\phi})\right].
$
Note that $Q=-|z^0|^2\phi$ is the $\gt$-length squared of the $(1,0)$ Euler
field $z^0\pa_0$.
\medskip

Define differential operators in a neighborhood of $M$ in $\C^{n+1}$ by 
$$
\Delta_{w,w'}=
\phi^{-1}\Delta_g + \frac{w'}{1-r\phi}\xi^a\pa_a
+ \frac{w}{1-r\phi}\xi^{\bb}\pa_{\bb} - \frac{ww'r}{1-r\phi}.
$$
By the observation in the proof of Proposition~\ref{laplprop},
$\Delta_{w,w'}$ has coefficients smooth near $M$.
It follows from
Proposition~\ref{laplprop} that extending a section of $\E(w,w')$ to be
harmonic with respect to $\Dt$ is the same as extending the corresponding
function on $M$ to be annihilated by $\Delta_{w,w'}$.  Therefore we
conclude:
\begin{proposition}\label{obschar}
If $w+w'+n+1=k \leq n+1$, the invariant operator $P_{w,w'}$
is the (normalized) obstruction to extending
a smooth function on $M$ to a smooth solution of
$\Delta_{w,w'}u=0$, where $\phi$ is chosen to solve
$J(\phi)=1+O(\phi^{n+2})$. 
\end{proposition}
\noindent
In particular, $P_{0,0}$ is the obstruction to smoothly solving $\Delta_g u
=0$ with $u$ having prescribed boundary value.  

In the setting of Proposition~\ref{obschar}, the obstruction, or
compatibility operators, were studied  
in \cite{g1}, \cite{g}, \cite{GL}.  In \cite{g}, the operators $P_{w,w'}$ 
were calculated explicitly for the sphere and the Heisenberg group with
their usual defining functions, which solve $J(u)=1$ exactly.  Of 
course in this flat situation, the ambient metric is invariantly defined to
all orders, so the invariant operators $P_{w,w'}$ exist without the
restriction $k\leq n+1$.
The result is that each $P_{w,w'}$ is a product of various of
the Folland-Stein operators $\Delta_b+i\al T$ on the Heisenberg group, or
the analogous operators of Geller on the sphere.  One can explicitly check
that for the Heisenberg group, the operator $P_{w,w'}$ in \cite{g} agrees
with the operator $\P_{w,w'}=(-2\Box)^k$ derived  
via tractors for flat pseudohermitian structures in
Proposition~\ref{flatgoody}. 

In \cite{GL}, a connection on a neighborhood of $M$ in $\C^{n+1}$
associated to an arbitrary defining function $\phi$ was
derived along with an explicit formula for $\Delta_g$ in terms of this
connection (Proposition 2.1 of \cite{GL}).
This together with the expressions for the curvature of this connection
given in \cite{GL} (for commuting derivatives)
in principle enable explicit calculation of the obstruction operators 
in the curved case in terms of pseudohermitian curvature and torsion and
their derivatives, and the transverse curvature $r$ and its derivatives.
For general $\phi$, the obstruction operator was calculated in \cite{GL} 
for $n=1$ and $w=w'=0$  to be:
\begin{equation}\label{P00}
\Delta_b^2 + T^2 +4\Im \nd_{\be}(A^{\al\be}\nd_{\al}),
\end{equation}
where all the quantities involved are those for the pseudohermitian
structure induced by $\phi$ on $M$. The CR invariance of this
operator was observed in \cite{hirachi}, from which one concludes
that for $n=1$, $P_{0,0}$ is given by this formula for any choice of
pseudohermitian structure.  
Direct calculation shows that the operator $\P_{0,0}:\ce(0,0)\to\ce 
(-2,-2)$ discussed after the proof of Proposition~\ref{flatgoody} agrees
with $P_{0,0}$. 

When $w=w'$, the invariant operators
can also be characterized as obstructions, or equivalently as log term
coefficients, when solving 
$(\Delta_g +w(n+1+w))u=0$.  This is analogous to the corresponding
characterization of the conformal operators $P_k$ given in
\cite{GZ}.  The resolvent $(\Delta_g - s(n+1-s))^{-1}$ has been studied 
in \cite{EMM}.
\begin{proposition}
If $u$ is a function on $\C^{n+1}$, then
\begin{equation}\label{scat}
\Dt(|z^0|^{2w}\phi^w u )
= (|z^0|^{2}\phi)^{w-1}(\Delta_g + w(n+1+w))u.
\end{equation}
Therefore, if $n+1+2w\in \N$ and $w\leq 0$, then $P_{w,w}f$ is the
(normalized) obstruction to solving $(\Delta_g + w(n+1+w))u = 0$ with
$\phi^wu$ smooth and $\phi^wu= f$ on $M$, where $\phi$ is chosen to solve 
$J(\phi)=1+O(\phi^{n+2})$. 
\end{proposition}
\begin{proof}
The relation \eqref{scat} is a direct calculation from \eqref{homoglapl}.
One expands the right hand side of \eqref{homoglapl} applied to $\phi^wu$
using the Leibnitz
rule and simplifies using \eqref{def}, \eqref{r}, \eqref{raise},
and the fact that
$$
g^{a\bb}\phi_{a\bb} = \phi(1-r\phi)^{-1}((n+1)r\phi -n).
$$
This latter is a consequence of
substituting \eqref{ke} into $g^{a\bb}g_{a\bb}=n+1$.  
\end{proof}
\noindent
We remark that there is a similar relation in the case $w\neq w'$ 
in terms of a
family of modifications of $\Delta_g$ by both first and zeroth order terms,
obtained by calculating $\Dt((z^0)^w (\zb^0)^{w'}\phi^{(w+w')/2} u )$. 

We close with a brief discussion of CR $Q$-curvature, 
considered also in \cite{FH}.  First recall
Branson's (\cite{Br2}) formulation of 
$Q$-curvature in conformal geometry. For this discussion 
denote by $P^N_k$ the operator $P_k$ in dimension $N$.  Fix $k \in \N$.  
The construction of \cite{GJMS} shows that the operator $P^N_k$
is natural in the strong sense that $P^N_k f$ may be written as a linear
combination of complete contractions of products of covariant derivatives 
of the curvature tensor of a representative for the conformal structure
with covariant derivatives of $f$, with coefficients which are rational in
the dimension $N$.  Also, the zeroth order term of $P^N_k$ may
be written as $P^N_k1= (N/2-k) Q^N_k$ for a scalar Riemannian invariant
$Q^N_k$  
with coefficients which are rational in $N$ and regular at $N=2k$.  The
$Q$-curvature in even dimension $N$ is then defined as $Q=Q^N_{N/2}$.  
An analytic continuation argument in the dimension then shows that under
the conformal rescaling $\hat{g}=e^{2\Up}g$, we have 
$e^{N\Up}\hat{Q}=Q + P_{N/2}\Up$.

If $M$ is a pseudohermitian manifold, the $Q$-curvature of its Fefferman
metric is invariant under rotations, so defines a function on $M$ which we
denote by $Q_F$.  It is then natural to define the CR $Q$-curvature of 
$(M,\theta)$ by $Q_{\theta} = 2^{-n}Q_F$.  Observations above imply
that $Q_{\theta}$ is given by a universal formula in the pseudohermitian
metric, torsion, curvature and their covariant derivatives.  The
transformation law for conformal $Q$-curvature gives immediately that if
$\hat{\theta} =  e^{\Up}\theta$, then 
\begin{equation}\label{CRQ}
e^{(n+1)\Up}Q_{\hat{\theta}} = Q_{\theta} + P_{0,0}\Up.
\end{equation}
For $N=4$, the conformal $Q$-curvature is given by
$6Q=\Delta K +K^2 -3|\Ric|^2$, where $K$ denotes the scalar curvature.  
Using Lee's formulae for $\Delta$, $K$, and $\Ric$, it is straightforward
to calculate this for the Fefferman metric; one obtains that the
$Q$-curvature for 3-dimensional CR manifolds is given by:
\begin{equation}\label{3dQ}
3Q_{\theta} = 2(\Delta_b R -2\Im \nabla^{\al}\nabla^{\beta}A_{\al\beta}).
\end{equation}
This quantity was introduced by Hirachi \cite{hirachi}, who showed that it
gives the coefficient of the $\log$ term in the Szeg\"o kernel.  He also  
established the transformation law \eqref{CRQ} in terms of the operator
\eqref{P00} by direct calculation.  

The role of $Q$-curvature in CR geometry is not clear.
In the conformal case, the total $Q$-curvature $\int Q$
is an interesting invariant of a compact conformal 
manifold.  However, since the expression
in \eqref{3dQ} is a divergence, it follows that we have $\int_M
Q_{\theta}=0$ for the corresponding integral for 3-dimensional CR
manifolds.  Also, the $Q$-curvature vanishes identically for a hypersurface
in $\C^{n+1}$ with pseudohermitian structure induced by a solution of 
$J(\phi)=1 + O(\phi)$.  A proof of this fact using a 
characterization of $Q$-curvature in terms of the ambient metric is given 
in \cite{FH}.  It can also be seen directly from the definition by noting
that the zeroth order 
term $P_k^N1$ for the conformal operator for the Fefferman metric is a
multiple of $({\tilde{\Delta}}^k|z^0|^{k-N/2})|_{z^0=1}$.  However, already 
${\tilde{\Delta}}|z^0|^{k-N/2}=-{\tilde g}^{0{\ol 0}}\pa^2_{0{\ol
0}}|z^0|^{k-N/2}$ vanishes to second order at $k=N/2$.

\end{document}